\documentstyle[12pt,subeqnarray,amssymb]{article}%

\newcommand{\al}{\alpha }
\def\b{\beta }
\newcommand{\g}{\gamma }        
\def\d{\delta }             

\newcommand{\et}{\eta }
\newcommand{\th}{\theta }       
\newcommand{\k}{\kappa }
\def\l{\lambda }            
\newcommand{\m}{\mu }
\newcommand{\n}{\nu }
\newcommand{\x}{\xi }         
\newcommand{\p}{\pi }         
\newcommand{\r}{\rho }
\newcommand{\s}{\sigma }        
\def\t{\tau }
\newcommand{\f}{\phi }         
\def\ph{\varphi }
\newcommand{\ch}{\chi }
\newcommand{\ps}{\psi }         
\newcommand{\om}{\omega }       \newcommand{\Om}{\Omega }

\def\R{\Bbb R }     
\def\C{\Bbb C }
\def\Z{\Bbb Z }
\def\1{\bold 1 }
\def\K{\Bbb K }
\def\H{\Bbb H }
 \def\vs #1{\vspace*{#1cm} }         \def\hs #1{\hspace*{#1cm} }
      \newcommand{\ovl}{\overline}
\newcommand{\frs}{\frenchspacing\ }

 \def\lb #1{\label{#1} }
 \def\slb #1{\slabel{#1} }
 \def\rep #1{(\ref{#1})}
\newcommand{\bde}{\begin{description} }
\newcommand{\ede}{\end{description}}
\newcommand{\bib}{\bibitem}
\newcommand{\bit}{\begin{itemize} }
\newcommand{\eit}{\end{itemize}}
\newcommand{\bea}{\begin{eqnarray}}   \newcommand{\bean}{\begin{eqnarray*}}
\newcommand{\eea}{\end{eqnarray}}     \newcommand{\eean}{\end{eqnarray*}}
\newcommand{\bse}{\begin{subeqnarray}} 
\newcommand{\ese}{\end{subeqnarray}} 
\newcommand{\non}{\nonumber}
\def\||{\Vert }
 \def\inv #1{#1^{-1}}
 \def\dtd #1{\delta _{#1} }        
 
\def\cop{\triangle }
\newcommand{\fal}{\forall }
\newcommand{\apli}{\rightarrow }    \newcommand{\map}{\mapsto }
\newcommand{\imp}{\Rightarrow }
\newcommand{\bra}{\langle }        \newcommand{\ket}{\rangle }
\def\pd{\otimes }             \newcommand{\pds}{\oslash }  
\def\+{\oplus }
\newcommand{\pex}{\wedge }
 \def\cali #1{{\cal #1} }
\def\Tr{\mbox{Tr} }
\def\Rm{{\R}_{+} }
\def\RmR{\Rm\pds{\R} }
\def\G{G_{hp} }
\def\lG{{\cal G}_{hp} }
\def\iRm #1{\int_{\Rm}\frac{d#1}{#1}\, }
 
\def\Cal{C$^*$-algebra }      
\def\vN{von Neumann }
\def\Hvn{Hopf-von Neumann }
\def\BH{{\cal B}(H) }
 \def\del #1{\partial_{#1} }

\begin{document}

\title{\bf Fourier Duality as a Quantization\\
    Principle\thanks{Work supported by CAPES and CNPq, Brazil.}}
\author{R. Aldrovandi\thanks{e-mail: ra@axp.ift.unesp.br}\ \ and L.A.
Saeger\thanks{e-mail: saeger@qcd.th.u-psud.fr}\\ Instituto de F\'\i sica
Te\'orica - UNESP\\ Rua Pamplona 145\\ 01405-900 - S\~ao Paulo, SP\\ Brazil}
\date{}
\maketitle

{\small The Weyl-Wigner prescription for quantization on Euclidean phase spaces makes
essential use of Fourier duality. The extension of this property to more
general phase spaces requires the use of Kac algebras, which provide the
necessary background for the implementation of Fourier duality on general
locally compact groups. Kac algebras -- and the duality they incorporate -- are consequently examined as candidates for a general quantization framework
extending the usual formalism. Using as a test case the simplest non-trivial
phase space, the half-plane, it is shown how the structures present in the
complete-plane case must be modified. Traces, for example, must be replaced by their noncommutative generalizations -- weights -- and the correspondence
embodied in the Weyl-Wigner formalism is no more complete. Provided the
underlying algebraic structure is suitably adapted to each case, Fourier
duality is shown to be indeed a very powerful guide to the quantization of
general physical systems.}

\vs{.2}


\vs{.2}

Internat. J. Theor. Phys. {\bf 36}(2), 345-383 (1997) \hspace*{\fill}funct-an/9608005

\newpage

\section{Introduction}

The complete quantum description of a physical system presupposes the
identification of the space of observables, the scene of its dynamical
evolution. In the classical description, this space corresponds to a
subalgebra of the algebra $C^{\infty}(M)$ of infinitely differentiable
functions on the phase space. When the latter is the linear space ${\R}^2$ (or
${\R}^{2n}$), the relationship between the quantum and the classical cases is
well known, given as it is by the Weyl correspondence prescription. The
correspondence assumes the existence (and knowledge) of a possible classical
version of the system and of a general prescription -- ``quantization'' -- to
transform the classical into the quantum description. The ultimate goal is to
uncover some ``grand principle'', a rule providing directly the quantal
description: given a system, we should be able to identify the observables
and their space without the mediation of classical quantities. It would then
be possible to describe even purely quantum systems, for which there are no
classical limits. Such an objective is still far ahead and for the present
time we are condemned to proceed from classical systems, trying to work up a
general procedure of quantization from particular examples. This
trial-and-error approach blends rigorous assumptions with inferences from
previous case-study experience.

The Weyl prescription, in its original form \cite{wey}, makes use of a very
particular kind of duality, the Pontryagin duality which holds only when the
underlying group of linear symplectomorphisms is Abelian. The Pontryagin dual
of an Abelian group $G$ is the space of characters, which is also an Abelian
group, though not necessarily the same. Thus, the Euclidean spaces ${\R}^n$ and
the cyclic groups ${{\Z}}_n$ are self-dual in this sense, but the circle and the
group of integers are dual to each other. The link between these groups is
provided by generalized two-way Fourier transforms mapping the $L^1$-space of
one into the $L^{\infty}$-space of the other, and for this reason we shall
use the expression {\em Fourier duality} as a synonym to (eventually
generalized) Pontryagin duality. In the Euclidean linear case, the Abelian
group involved is formed by the translations on phase space, which is
isomorphic to its own Fourier dual: the Fourier transforms of functions on
${\R}^2$ are functions on ${\R}^2$. Given the Fourier transform $\tilde f$ of a
classical dynamical function $f$, the Weyl prescription yields the
corresponding operator (q-number function) as a Fourier transform of $\tilde
f$ with a projective operator kernel. A formal inverse procedure (first
considered by Wigner \cite{wig}), involving an integration on an operator
space, gives then the corresponding c-number function \cite{hil}. These
c-number functions and their Fourier transforms belong to twisted
(noncommutative) algebras, different from the usual Abelian algebras of
convolution and pointwise product. The twisted convolution on $L^1({\R}^2)$ and
the Moyal product on $L^{\infty}({\R}{^2})$ arise naturally from the projective
operator product through the Weyl correspondence. Because of this ``quantum''
origin, the deformation of the usual algebras of functions on phase space
they represent is considered a quantization \cite{bffls}.

The main difficulty comes from the fact that, for general systems, including
those whose phase space is the Euclidean space, the group acting on phase
space (a group of linear symplectomorphisms, here called {\em special
canonical group}) is not Abelian and/or compact. On compact groups the
integration implied in the Fourier transform is defined in a simple way, as
there exists a unique Haar measure, which is both left- and right-invariant.
Amongst non-compact groups, the existence of Haar measures is assured only
for those which are locally compact, though in general the left-invariant and
the right-invariant measures differ (when they happen to be equal, the group
is said to be {\em unimodular}). Group-to-group duality is, however,
restricted to the commutative case: the space dual to a non-Abelian group is
no more a group, but an algebra. Duality must be understood no more as a
relationship between groups, but as a relationship in a wider category. A
fair formulation for the general locally compact case has been obtained in
the early seventies, and led to the introduction of Kac algebras. These are
Hopf-von Neumann algebras with peculiar generalized measures, called Haar
weights. Actually, for non-unimodular locally compact groups, Fourier
duality is only possible in the Kac algebra framework. We must abstract from
groups to Kac algebras in order to have a Fourier duality. In this sense the
usual Weyl correspondence is part of a highly nontrivial {\em
projective duality} for the Abelian group ${\R}{^2}$, where an algebra generated
no more by linear, but by projective operators comes into play \cite{ralas2}.

Our objective here is to give a step further in the question of quantization
through the study of these analytic-algebraic aspects. The algebraic facet is
better known: it is necessary to resort to Hopf-von Neumann algebras. These
algebras are, however, rather involved operator algebras, on which many
different topologies and measures can be defined. The analytic facet lies
precisely in the choice of the correct topology and measure. Our guiding idea
will be the assumption of Fourier duality, which stands at the heart both of
the Weyl quantization approach and of the group duality alluded to above.
Since Fourier duality in its more general form is implemented in the Kac
algebra structural frame \cite{ensc}, we argue that they are also able to
provide a generalized Weyl prescription for quantizing a phase space on which
a separable locally compact type~I group acts by symplectomorphisms.

We take as a test case the simplest non-trivial example of phase space: the
half-plane. The fact that the configuration space $\Rm$ is Abelian, that the
manifold of the special canonical group coincides with the phase space
manifold and that there is no need to consider central extensions of this
canonical group accounts for the relative simplicity of the example. Its
non-triviality comes from the non-trivial properties of the group which,
besides being non-Abelian and non-compact, is non-unimodular. These are
important features, which bring to light the main difficulties of
quantization on a general phase space. Specifically, this example also shows
why the usual Weyl-Wigner quantization procedure does not
generalize straightforwardly and does not always lead to a generalized Moyal
bracket \cite{moy}, or to a deformation of the algebra of functions on phase
space. Although this example does not cover quantization on a general phase
space, where the respective canonical group may have little to do with the
space it acts \cite{ish}, we believe the duality principle it illustrates
can be generalized to quantization on any phase space to which an operator
algebra can be associated, as it is done in \cite{land}.

We begin with an exposition of the classical picture on the half-plane in
section~\ref{canon}. We use, in order to select a group on phase space, Isham's
canonical approach which, though not quite general, is enough for the case in
view. In the next section we give some details and classify the induced
irreducible representations of the half-plane special canonical group, which is
in fact a special parametrization of the affine group on the line (conversely,
we show in the Appendix that the half-plane is the unique non-trivial
homogeneous symplectic manifold of the affine group on the line). Since it
seems that neither Hopf-von Neumann nor Kac algebras are structures quite
familiar to the Physics community, we review them in a separate section. We
emphasize those Kac algebras which are related to groups, in order to show how
group duality is attained. At the end of that section we also show how to
decompose the operator Kac algebra of a type~I group according to its unitary
dual. This is not found in the Kac algebra literature and will be essential to
the interpretation of the Weyl formula in the duality framework. The half-plane
case, used all along more as a gate into non-trivial aspects, is finally
retaken for its own sake and given its finish in section~\ref{qohpl}. The whole
treatment leads to a reappraisal of the reach and limitations of the
Weyl-Wigner formalism as a guide for quantization on general phase spaces.

\section{The Special Canonical Group}
\lb{canon}

The phase space we want to quantize on is the half-plane $\Rm\times{\R}$ the
cotangent bundle of the configuration space given by the half-line $\Rm$. On
this symplectic manifold we use the coordinates $x$ and $p$, in terms of
which the symplectic form, given as the derivative of the Liouville canonical
1-form, is
\bean
 \om = d\th_o=dp\pex dx,\hs 1x\in\Rm,\, p\in{\R}.
\eean
The symplectic form implements, through the equation
\bea
 i_{X_f}\om=-df,
\lb{eqm}
\eea 
a homomorphism between the space of $C^{\infty}$-functions (Hamiltonians) and
the space of symplectic Hamiltonian vector fields, whose kernel are the
constants. Since $\om$ is non-degenerate, \rep{eqm} can be solved for the
vector fields and yields, in the above coordinates, $X_f=\del pf\,\del x-\del
xf\,\del p$. The symplectic form also provides a Lie algebra structure on
the $C^{\infty}$-functions, as it defines the Poisson bracket by
\bea
 \{f,g\}=-\om(X_f,X_g),
\lb{dpbr}
\eea
which is isomorphic to the Hamiltonian vector fields Lie algebra through
\bean
[X_f,X_g]=-X_{\{f,g\}}.
\eean 

We shall follow Isham \cite{ish} in the first steps. To quantize a phase space
we start by looking for a finite dimensional (for simplicity) group whose
elements act as symplectomorphisms, that is, preserve the symplectic structure.
The action must be transitive, so as to avoid any lack of globality in the
quantum description, and also (quasi-)effective. It is thus necessary to find a
finite dimensional group $\G$ under whose action the half-plane is a symplectic
homogeneous G-space. The task can be simplified by proceeding as follows.
Consider a group whose manifold is the configuration space, $(\Rm,\cdot)$, and
make it act on a linear space so as to get at least an almost-faithful
representation (a representation whose kernel is discrete). Take the action of
$\Rm$ on ${\R}$ given by
\bea
 \l\in\Rm\map R_{\l}(a)=a\l,\ a\in{\R}.
\lb{repl} 
\eea
Construct the semi-direct product group $\RmR$, with the product operation
given by
\bean
 (\l,a)(\r,b)=(\l\r,a+\f_{\l}(b)),
\eean
where $\f_{\l}(b)=R_{\l^{-1}}^*(b)=b/\l$ is the homomorphism on ${\R}$ given by
the representation $R^*$ contragradient to \rep{repl}. The identity in
$\RmR$ is $(1,0)$ and the inverse element of $(\l,a)$ is given by
$(\l^{-1},\f_{\l}^{-1}(-a))$, with $\f_{\l^{-1}}=\f_{\l}^{-1}$.

Considering the left action
\bean
l_{(\l,a)}(x,p)=(\l x,p/\l - a),
\eean
of this group on the space $\Rm\times{\R}$, we see that $\G$ = $\RmR$ is 
formed by some special linear canonical transformations on the half-plane.
Actually, we show in the Appendix, using Kirillov's orbits method \cite{kir},
that the half-plane is the only non-trivial symplectic manifold canonically
invariant by $\G$.

The Lie algebra $\lG$ of $\G$ can be obtained from the group product with the
help of the formula $e^{tA}e^{sB}e^{-tA}e^{-sB}=e^{ts[A,B]+\ higher\ orders\
in\ t,s}$ and is given on ${\R}\+{\R}$ by
\bea
 [(l,a),(r,b)]=(0,ar-lb).
\lb{lab}
\eea
It is straightforward to realize this Lie algebra in terms of symplectic
Hamiltonian vector fields on the half-plane. By the exponential mapping
$(l,a)\map (e^l,a)\in\RmR$ we introduce the 1-parameter subgroups $t\map
(e^{lt},0),\ s\map(1,as)$, whose action on $\Rm\times{\R}$,
\bean
 l_{l,a}(x,p)=(e^{lt}x,e^{-lt}p-as)
\eean
is easily found to be generated by the symplectic Hamiltonian
(right-invariant) vector fields
\bean
 X_{l,a}(x,p)=lx\del x -(lp+a)\del p,
\eean
corresponding to the Hamiltonians $h_{l,a}(x,p)=ax+lxp$. On
$C^{\infty}(\Rm\times{\R})$ these Hamiltonians define a Poisson subalgebra by
\bea
\{h_{l,a},h_{r,b}\}=h_{0,ar-lb},
\lb{pbhf}
\eea 
whose structure is identical to that of \rep{lab}. We can then say that there
is a faithful momentum mapping $J:T^*\Rm\apli\lG^*$, allowing the association
of the pair $(l,a)\in\lG$ to the Hamiltonian function $h_{l,a}$ by the
duality pairing $\bra J(x,p),(l,a)\ket=h_{l,a}(x,p)$. By this Lie algebra
isomorphism we privilege the functions $h_{l,a}$ as a preferred class of
observables to be quantized. Also because of this isomorphism, there is no
need to central-extend $\lG$ as it happens in the complete-plane case. In
other words, the cohomology space $H^2(\lG,{\R})\sim H^2(\G,{\R})$ is trivial
\cite{ish}. It is then possible to take the unitary irreducible linear
representations of the special canonical group realized in terms of operators
on a given Hilbert space and try to find an unbounded operator
(representation generator) in correspondence with each preferred observable
on the half-plane. In the next section we will provide the representations
necessary to characterize such operators but, differently from Isham's
approach, functions will be associated to bounded operators {\em a la} Weyl,
which means that we shall work at the group representation level.

\section{Induced Representations and the Unitary Dual}
\lb{irrur}

Irreducible unitary representations of semi-direct product groups are easily
constructed via Mackey's induced representation theory (see
\cite{mack,bar,gur,sugi}). In this section we construct irreducible unitary
representations of $\G=\Rm\pds{\R}$ by that method. We first note that a $\G$
element $g=(\l,a)$ can be decomposed in its $\Rm$ and ${\R}$ parts according to
\bea
(\l,a)=(\l,0)(1,\f_{\l}^{-1}(a))=g_{\Rm}g_{\R}
\lb{dlaie}
\eea
and 
\bea
(\l,a)=(1,a)(\l,0).
\lb{dlaid}
\eea

We begin by looking for unitary irreducible representations of the subgroup
${\R}$. This is immediate, since ${\R}$ is an Abelian normal subgroup. Its
character (one-dimensional) representations on ${\C}$ are given by
$V_x(a)=e^{ixa}$, where $x$ is a label contained in the unitary dual group
$\hat{\R}\sim{\R}$.
The Hilbert space where our induced representation of $\G$ will be realized
is constructed as the space of functions $f:\G\apli{\C}$ which can be
decomposed into wavefunctions $\x$, $f_x(g)=V_x^{-1}(g_{\R})\x(g_{\Rm})$,
or, using \rep{dlaie} with $g=(\l,a)$, 
\bea
f_x(\l,a)=e^{-ix\f_{\l}^{-1}(a)}\x(\l),\hs 1\x\in L^2(\Rm),
\lb{hsd}
\eea 
and on which $f_x$ satisfies
\bean
 \iRm{\l}\ |f_x(\l,a)|^2=\iRm{\l}\ |\x(\l)|^2<\infty.
\eean
We indicate this space by $H_x(\G)$ and, in agreement with \rep{hsd}, use the
fact that it is isomorphic to $L^2(\Rm)$. The induced representation of $\G$ on
$H_x(\G)$ is then defined, for each $x\in{\R}$, by $[T_x(g)f_x](h)=f_x(\inv gh)$,
or, directly in terms of wavefunctions in the coordinate representation (that
is, on $L^2(\Rm)$), by
\bea
[T_x(g)\x](h_{\Rm})=V_x^{-1}([g^{-1}h_{\Rm}]_{\R})\x([g^{-1}h_{\Rm}]_{\Rm}),
\lb{rin}
\eea 
or still, even more explicitly, using $g=(\l,a),\ h=(\r,b)$, and computing
$g^{-1}h_{\Rm}=\linebreak(\l^{-1}\r,\f_{\r}^{-1}(-a))$,
\bea
[T_x(\l,a)\x](\r)=e^{ix\f_{\r}^{-1}(a)}\x(\l^{-1}\r).
\lb{uiir}
\eea
That these operators do represent the group $\G$, 
\bea
 T_x(\l,a)T_x(\r,b)=T_x((\l,a)(\r,b)),
\lb{lrr}
\eea 
follows trivially from comparing the two identities below: applying the
left hand side of \rep{lrr} to $\x\in L^2(\Rm)$, we obtain
\bean
[T_x(\l,a)T_x(\r,b)\x](\et)&=&e^{ix\f_{\et}^{-1}(a)}[T_x(\r,b)\x](\l^{-1}\et)\\
&=&e^{ix\f_{\et}^{-1}[a+\f_{\l}(b)]}\x((\l\r)^{-1}\et),
\eean
while the right hand side gives
\bean
[T_x((\l,a)(\r,b))\x](\et)&=&[T_x(\l\r,a+\f_{\l}(b))\x](\et)\\
&=&e^{ix\f_{\et}^{-1}[a+\f_{\l}(b)]}\x((\l\r)^{-1}\et).
\eean
Unitarity and irreducibility of the representation \rep{uiir} will be proved in
the following.

Abstracting from the Hilbert space $L^2(\Rm)$, we can write the operator
$T_x$ as
\bean
 T_x(\l,a)|_{\r}=e^{ix\f_{\hat\r}^{-1}(a)}e^{-i\ln(\l)\hat\p},
\eean
where $|_{\r}$ means that the operator acts on the argument $\r$ in such a
way that the multiplication and dilation operators are defined by
\bean
\hat\r \x(\r)&=&\r \x(\r)\\ 
\hat\p \x(\r)&=&-i\r\del{\r}\x(\r).
\eean

An operatorial version of the decompositions given at the beginning of this
section comes up if we define the operators (dropping $|_{\r}$ from now on)
\bse
T_x(\l,0)&\equiv& L(\l)=e^{-i\ln(\l)\hat\p};\slb{opsa}\\ 
T_x(1,a)&\equiv&\hat V_x(a)=e^{ix\f_{\hat\r}^{-1}(a)}\slb{opsb},
\lb{ops}
\ese
with which we have
\bea
 T_x(\l,a)=\hat V_x(a)L(\l),
\lb{eto} 
\eea
where $L(\l)$ is identified as the left-regular unitary representation of the
group $(\Rm,\cdot)$ on $L^2(\Rm)$. Definitions \rep{ops} also allows us to
rewrite \rep{dlaie} and \rep{dlaid} in operatorial form
\bea
 L(\l)\hat V_x(\f_{\l}^{-1}(a))=\hat V_x(a)L(\l).
\lb{rcl}
\eea 
Expanding the identity above according to \rep{ops}, and recalling that
$\f_{\l}(a)=a/\l$, we obtain, up to first order in $al,\ l=\ln\l$,
\bean
 [\hat\r,\hat\p]=i\hat\r.
\eean

Now the unitarity of \rep{uiir} follows easily from \rep{eto} and the
unitarity of $L$ and $V_x$,
\bean
 T_x^{\dagger}(\l,a)&=&L(\l^{-1})\hat V_x(-a)\\
 &=&\hat V_x(\f_{\l}^{-1}(-a))L(\l^{-1})\\ 
 &=&T_x((\l,a)^{-1}),
\eean
where the second equality comes from \rep{rcl}.

At this point we should ask whether there exists an equivalence relation
between the operators $T_x$. This is an important question if we want to do
harmonic analysis on the group $\G=\Rm\pds{\R}$, as we shall, for we must sum
(integrate) over the unitary dual $\widehat{\G}$ of $\G$, the space of classes
of inequivalent irreducible representations. To answer the question we begin
by observing that the right-regular representation of $(\Rm,\cdot)$ acts on
$\x\in L^2(\Rm)$ by $[R(\r)\x](\et)=\x(\et\r)$. In order to verify whether
this operator is an intertwining for the $T_x$, we calculate
\bea
[R(\r)^{-1}T_x(\l,a)R(\r)\x](\et)&=&[T_x(\l,a)R(\r)\x](\et\r^{-1})\non\\
&=&e^{ix\f^{-1}_{\et\r^{-1}}(a)} \x(\l^{-1}\et).
\lb{intw}
\eea
Now, remembering that $\Rm$ acts on ${\R}$ by $\f_{\r}$, its associated
co-action $\hat{\f}_{\r}$ on the character space $\hat{\R}$ is
defined by
\bean
[\hat{\f}_{\r}\ch_x](a)&\equiv&\ch_x({\f}_{\r}^{-1}(a))=e^{ix\f_{\r}^{-1}(a)}.
\eean
With this at hand, it is easy to see that the coefficient of the wavefunction
$\x$ in \rep{intw} above is just the following co-action:
\bean
[\hat{\f}_{\et\r^{-1}}\ch_x](a)=\ch_x({\f}_{\et\r^{-1}}^{-1}(a))&=&
\ch_x(\f_{\r}\circ\f_{\et}^{-1}(a))\\
&=&\hat{\f}_{\r^{-1}}\ch_x(\f_{\et}^{-1}(a)).
\eean
Explicitly, considering that $\ch_x\in\hat{\R}\sim{\R}\ni x$, we have the
co-action given by
\bea
 \hat{\f}_{\r^{-1}}(x)=\r^{-1}x.
\lb{coa}
\eea
We then conclude that the right-regular representation $R$ is an intertwining
operator for the $T$'s, connecting them by the co-action $\hat{\f}$,
\bea
R(\r)^{-1}T_x(\l,a)R(\r)=T_{\hat{\f}_{\r^{-1}}(x)}(\l,a).
\lb{clr}
\eea
We have then three classes of representations: one for $x>0$ and one for $x<0$,
both isomorphic to $\Rm$; and that one represented by the point $x=0$. We shall
indicate the cases $x>0$ or $x<0$ simply by $\pm$ and write the two
infinite-dimensional representation operators as
\bea
 T_{\pm}(\l,a)=e^{\pm ia\hat\r}L(\l).
\lb{idr}
\eea
In the case $x=0$, we have simply $L(\l)$, the left-regular unitary
representation of $(\Rm,\cdot)$, $T_0(\l,a)=L(\l)$. This representation is
reducible, that is, it is possible to decompose it in terms of the
$(\Rm,\cdot)$ characters $\ch_{y}(\l)=\l^{iy},\ y\in{\R}$, and write formally
\bean
L=\int_{\R}^{\+}dy\,\ch_{y}.
\eean
This give us an infinity $({\R})$ of 1-dimensional irreducible representations,
\bea
 T_{y}(\l,a)=\l^{iy}.
\lb{1dr}
\eea 
Summing up: once we suppose the irreducibility of the $T_x$, which will be
proved just below, the unitary dual $\widehat{\G}$ is given by
$\{+\}\cup\{-\}\cup{\R}$. If we compare this result with the orbits of the
coadjoint action of $\G$ obtained in the Appendix, we observe that the
formula $\widehat{\G}=\lG^*/\G$ holds.

Now, to address the problem of irreducibility of the induced representations
$T_x$ we shall refer to an important result of Mackey's theory. Mackey's
imprimitivity theorem \cite{bar,tay} for semi-direct products states that the
induced representations of such groups, in our case $T_x$, will be
irreducible if and only if the semi-direct product group $\Rm\pds{\R}$ which it
represents satisfies a condition of {\em regularity}. This condition
essentially means that the $\Rm$-orbits in $\hat{\R}$ by the $\hat{\f}$
action are countably separated with respect to the Borel structure. This is
easily seen to be fulfilled since $\hat{\R}={\R}_-\cup\{0\}\cup\Rm$. So, our
group $\G$ is regular and the representations $T_x$ are irreducible.

The above analysis gives still another important information about the group
$\G$. Type~I groups are, roughly speaking, those groups which have a well
behaved Borel structure on the unitary dual, more specifically, the
decomposition of representations of these groups into irreducible
representations is unique \cite{mack1}. Good examples are the Abelian and
the semi-simple groups. From another theorem by Mackey \cite[p. 536]{bar}, a
regular semi-direct product group, say $\Rm\pds{\R}$, is a {\em type~I group}
if and only if for each $x\in\hat{\R}$, its isotropy subgroup $I_x$ is a
type~I group. Well, we know that the orbits through $x$ are given by ${\cal
O}_x=\Rm/I_x$, and we have found that they are isomorphic either to $\Rm$
$(\Rm,\cdot)$ or to the trivial $\{e\}$. Consequently, each isotropy
subgroup is necessarily isomorphic to one of them, and they are both of type~I.

\section{Kac Algebras and Group Duality}

Once characterized and constructed the representations of the group under
which the half-plane is canonically invariant, we must give a rule to
associate an operator to each observable. To do it we will use the powerful
techniques provided by the Kac algebras. These algebras were constructed in
1973 independently by G.I. Kac and L.I. Va\u\i nermann \cite{kava}, and M.
Enock and J.-M. Schwartz \cite{ens}, with the objective of generalizing to
non-unimodular locally compact groups the Pontryagin (Abelian groups) and
Tannaka-Krein (compact groups) duality theorems. A duality for locally
compact (l.c. from now on) non-unimodular groups, comprising previous works
of P. Eymard, N. Tatsuuma and J. Ernst on a category wider than that of such
groups, had already been partially obtained in the seventies by M. Takesaki
\cite{tak} in the Hopf-von Neumann algebra framework. Unfortunately, due to
an incomplete theory of noncommutative integration, Takesaki's work on that
direction had a lack of symmetry. A general duality only was possible after
a considerable knowledge on weights was obtained. This knowledge led to the
definition of Kac algebras by the addition of a suitable (Haar) weight on \Hvn
algebras.

Actually, a general duality for locally compact groups is achieved if we
associate to them two Kac algebras, one on the von Neumann algebra of
$L^{\infty}$-functions and the other on the von Neumann algebra generated by
left-regular representations. These two algebras turn out to be dual in the
category of Kac algebras. This means that, by duality, to each
$L^{\infty}$-function on the group we can make to correspond an operator
written as a linear combination of the left-regular representations. In this
section we will introduce \Hvn and Kac algebras, apply the latter to groups in
order to show the l.c.\frs group duality, and show how they decompose following
the unitary dual of a type~I group.

Since Hopf-von Neumann and Kac algebras are, to begin with, von Neumann
algebras, we start by recalling some definitions on these algebras which will
be necessary (see, for example Ref.~\cite{bra}). A {\em von Neumann algebra}
$M$ is an {\em involutive unital} subalgebra of the Banach algebra $\BH$ of
bounded linear operators on a Hilbert space $H$, closed with respect to the
{\em strong topology}, a topology which is defined by the open balls of the
family of {\em seminorms}
\bea
 \|| a\||_{F,\ps}=\||a\ps\||,\ \ps\in H.
\eea
Besides this topology on $\BH$ there is also the {\em uniform} topology,
which is defined through the norm
\bea
 \Vert a\Vert=\sup\{\Vert a\ps\Vert_H,\ \Vert\ps\Vert_H=1,\ \ps\in H\}.
\lb{wt}
\eea 
As a map $\Vert\ \Vert :M\apli [0,\infty]$, this norm satisfies the following
conditions:
\bit
\item $\|| a\||=0$ if and only if $a=0$;
\item $\||{a+b}\||\leq\|| a\||+\|| b\||$;
\item $\||{\al a}\||=|\al|\|| a\||,\hs 1\al\in{\C}$;
\item $\Vert ab\Vert\le\Vert a\Vert\, \Vert b\Vert.$
\eit
The first axiom does not hold for a {\em seminorm}. We shall consider also an
another topology coming from seminorms, the {\em ultra}($\s$)-weak topology.
It is given by $\Vert T\Vert_{\s,\ps_i,\f_i}=\sum_i|(T\ps_i,\f_i)|$, where
$\ps_i,\f_i\in H$ are such that $\sum_i\||\ps_i\||^2<\infty$ and
$\sum_i\||\f_i\||^2<\infty$. The {\em predual} $M_*$ of $M$ is the (Banach)
space of the ultra-weakly continuous linear functionals on $M$.

The word {\em involutive}, used in the definition above, means that on $M$ is
defined a map $*:M\apli M$, the {\em involution}, such that
\bit
\item $(\al a+\b b)^*=\ovl{\al}a^*+\ovl{\b}b^*$;
\item $(ab)^*=b^*a^*$;
\item $(a^*)^*=a$.
\eit
Besides these axioms, on a von Neumann algebra it is also true that
$\||a^*\||=\|| a\||$ (it is an involutive Banach algebra), $\||a^*a\||=\||
a\||^2$ (it is a $C^*$-algebra) and that the unit is preserved by the
involution, ${\1}^*= {\1}$. Finally, a $W^*$ algebra is an algebra $M$ which equals
the dual of its predual, $M=(M_*)^*$. It is, roughly speaking, an abstract
$C^*$-algebra which can always be realized as a von Neumann algebra on a
suitable Hilbert space $H$.

We can introduce at this point the definition of a Hopf-von Neumann algebra.
A {\em co-involutive \Hvn algebra} is a triple ${\H}=(M,\cop,\k)$
where \cite{ensc}
\bit
\item $M$ is a $W^*$-algebra;

\item the homomorphism $\cop:M\apli M\pd M$, called {\em coproduct}, is
normal, injective and such that 
\bse 
\cop{\1}&=&{\1}\pd{\1}\\ 
(\cop\pd id)\circ\cop&=&(id\pd\cop)\circ\cop.\slb{cob}
\lb{coaa}
\ese
The first statement above says that $\cop$ is unital and the latter that it is
coassociative. Since $\cop$ is a homomorphism of $W^*$-algebras, this means
that it is linear and
\bea
\cop(a b)&=&\cop(a)\cop(b)\lb{hcm}.
\eea
\item there is a map $\k:M\apli M$, called {\em co-involution}, which is an
involutive anti-automorphism, that is, which is linear and such that, $\fal
a,b\in M$,
\bse
 \k(ab)&=& \k(b)\k(a);\slb{kanc}\\
 \k(a^*)&=&\k(a)^*;\slb{kanb}\\ 
 \k(\k(a))&=&a.\slb{kana}
\lb{kant}
\ese
\item it is also an anti-coautomorphism,
\bea
(\k\pd \k)\circ\cop =\s\circ\cop\circ \k,
\lb{acoa}
\eea
where $\s(a\pd b)=b\pd a$.
\eit

${\H}$ is said to be Abelian or commutative if M is Abelian, and symmetric or
cocommutative if $\s\circ\cop=\cop$. Note that from \rep{kana} and \rep{kanb}
it follows that $\k(\k(a^*)^*)=a$, or $\k\circ\ast\circ \k\circ\ast=id$, but
the converse is not true. This condition is actually weaker than those axioms.
One of the differences between Hopf-von Neumann algebras and Woronowicz's
``compact matrix pseudo-groups'' \cite{wor} is that this weaker condition is
imposed instead of \rep{kana}, \rep{kanb}.

Given a co-involutive Hopf-von Neumann algebra $(M,\cop,\k)$, where $M$
acts on the Hilbert space $H$, and a representation $\m$ of its predual $M_*$
on the Hilbert space $H_{\m}$, a partial isometry $U\in\cali B(H_{\m})\pd M$
such that
\bea
 \m(\om)=(id\pd\om)(U),\hs 1 \om\in M_*,
\lb{gor} 
\eea
is said to be the {\em generator} of $\m$. If $\m$ is {\em multiplicative}
and {\em involutive}, its generator $U$ satisfies the respective identities
\bse
 (id\pd\cop)(U)&=&(U\pd{\1})({\1}\pd\s)(U\pd{\1})({\1}\pd\s)\slb{apr}\\
 (id\pd\om\circ\k)(U)&=&(id\pd\om)(U^*).\slb{mein}
\ese
In the following we shall also denote $(U\pd{\1})$ and
$({\1}\pd\s)(U\pd{\1})({\1}\pd\s)$ in $\cali B(H_{\m})\pd M\pd M$ by $U_{12}$ and
$U_{13}$ respectively. If $\x,\et\in H$ we define the linear form
$\om_{\x,\et}\in M_*$ by
\bea
 \bra a,\om_{\x,\et}\ket\equiv (a\x|\et)_H,\hs 1\fal a\in M.
\lb{dlf}
\eea
The formula 
\bea
(\hat U(\al\pd\b)|\g\pd\d)_{H\pd H_{\m}}=(\b|\m(\om_{\g,\al})\d)_{H_{\m}},\hs
1\al,\g\in H,\ \b,\d\in H_{\m},
\lb{dfw}
\eea
coming from \rep{gor} and \rep{dlf}, and connecting the representation
$\m$ and the operator $\hat U\equiv\s\circ U^*\circ\s\in M\pd\cali B(H_{\m})$
(the dual of $U$), will be very useful.

Before introducing Kac algebras, some facts concerning weights and the
representation of a von Neumann algebra by a weight -- the GNS construction
-- are worth mentioning. The basic references are \cite{bra} and
\cite[section~2.1]{ensc}.

Consider a map from the set of strictly positive elements of
$M$, $\ph:M^+\apli [0,\infty]$, with the conditions
\bit
\item $\ph(a+b)=\ph(a)+\ph(b);$
\item $\ph(\l a)=\l\ph(a),\ \fal\l\ge 0,\ \mbox{where}\ 0\cdot\infty\equiv 0;$
\item $\ph(a^*a)=\ph(aa^*)\ \fal a\in M.$
\eit
The first two conditions define a {\em weight} on $M$, and the three together
define a {\em trace} \cite{dix2}. A weight generalizes the concept of positive
linear functional on C$^*$-algebras and, in particular, the concept of {\em
state}. Associated to $\ph$ we define the left ideal $\cali N_{\ph}\subset M$
by $\{a\in M\,|\,\ph(a^*a)<\infty\}$, and the involutive algebra $\cali
M_{\ph}$ as the linear span of $\{a\in M^+\, |\,
\ph(a)<\infty\}\subseteq\cali N_{\ph}\cap\cali N_{\ph}^*$, with $\cali
N_{\ph}^*=\{a^*\ |\ a\in \cali N_{\ph}\}$. A weight $\ph$ is called:
\bde
\item{\bf normal} if for every sequence $\{a_i\}$ with upper bound $a\in M^+$, 
$\ph(a)$ is the upper bound of the sequence $\{\ph(a_i)\}$;
\item{\bf faithful} if $\ph(a)=0\imp a=0,\ a\in M^+$,
\item{\bf semi-finite} if $\cali M_{\ph}$ is ultra-weakly dense in $M$.
\ede

Given a normal, faithful and semi-finite weight $\ph$ on a von Neumann
algebra $M$, we construct a representation of $M$ by the following procedure
\cite{bra}: $\ph$ defines a scalar product in $\cali N_{\ph}$, through
\bean
(a|b)_{\ph}\equiv\ph(b^*a).
\eean
It is actually only a quasi-scalar product since, as $\ph(a^*a)\geq 0,\
(a|a)_{\ph}$ can be zero. To circumvent this problem, we should factor the
left ideal $I_{\ph}=\{b\in A\ |\ (b|b)_{\ph}=0\}$ out of $\cali N_{\ph}$. The
quotient is formed by equivalence classes [b] of elements b' such that b-b'
is in $I_{\ph}$. $\cali N_{\ph}/I_{\ph}$ has a pre-Hilbert structure given by
the scalar product $ ([a],[b])=(a|b)_{\ph} $ which is invariant on each
class. Completing $\cali N_{\ph}/I_{\ph}$ with respect to this product we get
the Hilbert space $H_{\ph}$. The map
\bean
\p_{\ph}(a): \cali N_{\ph}/I_{\ph} & \apli & \cali N_{\ph}/I_{\ph}\\ \ [b] &
\map & [ab]
\eean
is bounded and can be extended to $H_{\ph}$ as a bounded operator. We call
($\p_{\ph},H_{\ph}$) the {\em GNS construction} of ($M,\ph$), and $a_{\ph}$
denotes the image of $a\in\cali N_{\ph}$ into $H_{\ph}$ by the canonical
injection $\p_{\ph}:\cali N_{\ph}\apli H_{\ph}$, $a\map\p_{\ph}(a)=[a]$. The
image $\p_{\ph}(\cali N_{\ph}\cap\cali N_{\ph}^*)$ is a {\em left Hilbert
algebra} \cite{ensc}, which is isomorphic to $M$. The image of the involution
$\ast$ is the operator $S_{\ph}$, which has the polar decomposition
\bean
 S_{\ph}=J_{\ph}\cop_{\ph}^{1/2}.
\eean
This decomposition gives rise to the antilinear isometry
$J_{\ph}:H_{\ph}\apli H_{\ph}$, such that $JMJ=M',\ JaJ=a^*,\ a\in Z(M)$, and
to the {\em modular operator} $\cop_{\ph}$, where $M'=\{a\in\BH\ |\ ab=ba,
for\ all \ b\in M\}$ is the {\em commutant} of $M$, and $Z(M)=M\cap M'$ is
the {\em center} of $M$. The modular operator satisfies
$\cop^{it}_{\ph}M\cop^{-it}_{\ph}=M$, for all $t\in{\R}$, and this leads to the
definition of the {\em modular automorphism group} $\s_t^{\ph}$ on $M$ by
\bea
\s_t^{\ph}(a)=\cop^{it}_{\ph}\,a\,\cop^{-it}_{\ph}.
\lb{modg}
\eea
The modular group $\s^{\ph}_t$ is such that the weight $\ph$ is invariant,
$\ph=\ph\circ\s_t^{\ph}$, and is also characterized by the fact that 
 $\ph$ is the unique KMS-weight associated to it. This short overview
on the Tomita-Takesaki theory extended to weights will be enough to introduce
Kac algebras.

A {\em Kac algebra} ${\K}=(M,\cop,\k,\ph)$ satisfies the following axioms (for
a good review on Kac algebras, see the first sections of an article by one of
its founders in Ref.~\cite{vai}):
\bit
\item ($M,\cop,\k$) is a co-involutive \Hvn algebra;
\item $\ph:M^+\apli [0,\infty]$ is a normal, faithful and semi-finite
weight on M called {\em Haar weight} such that
\bit
\item $\cop(\cali N_{\ph})\subset\cali N_{{\1}\pd\ph}$. A stronger version
of this axiom is more manageable and will be used. It says that $\ph$ is
left-invariant with respect to $\cop$, or
\bea
(id\pd\ph)\cop (a)=\ph(a){\1}\qquad \fal a\in M^+;
\lb{lia}
\eea
\item $\ph$ is symmetric, $\fal\ a,b\in\cali N_{\ph}$,
\bea
(id\pd\ph)[({\1}\pd b^*)\cop(a)]=
\k\circ(id\pd\ph)[\cop(b^*)({\1}\pd a)];
\lb{syma}
\eea
\item and
\bea
 \k\circ\s^{\ph}_t=\s^{\ph}_{-t}\circ \k\qquad \fal\ t\in{\R}.
\lb{mga}
\eea
\eit\eit

Given a \Hvn algebra ${\H}$, then a Haar weight which makes of ${\H}$ a Kac
algebra, if it exists, is unique up to a scalar \cite[sect. 2.7.7]{ensc}. A
Kac algebra ${\K}$ is called {\em unimodular} if the Haar weight $\ph$ is a
trace and is invariant by $\k$, $\ph =\ph\circ \k$. When $\ph$ is a
trace, then it is true that $\cop_{\ph}=1$ and $\s_t^{\ph}=id.$, as it
happens for the Abelian Kac algebras of groups described in the next
subsection.

Associated to a Kac algebra there exists always an isometry $W$ belonging to
$M\pd\cali B(H_{\ph})$, called the {\em fundamental operator}, such that
\bea
 W(a_{\ph}\pd b_{\ph})=[\cop(b)(a\pd{\1})]_{\ph}\qquad a,b\in\cali N_{\ph}.
\lb{cpw}
\eea
This unitary operator implements the coproduct as follows:
\bea
 \cop(a)=W({\1}\pd a)W^*.
\lb{canimcop}
\eea

Let us now introduce the dual of a Kac algebra ${\K}$ based on $M$. Its predual
$M_*$ has a product $\ast$ given by
\bea
 \bra a,\om\ast \om'\ket=\bra\cop a, \om\pd \om'\ket,
\lb{dpcp}
\eea
and an involution $^o$ by
\bea
 \bra a,\om^o\ket=\ovl{\bra\k(a)^*,\om\ket}.
\lb{rei}
\eea
The predual is thus an involutive Banach algebra. Besides the GNS
representation of $M$ on $H_{\ph}$, there is a multiplicative and involutive
representation of $M_*$ on the same Hilbert space,
\bean
 \l:M_*\apli\hat M\subset\cali B(H_{\ph}).
\eean
The representation $\l$ is such that $\l(\om)$ is a bounded operator on
$\p_{\ph}(\cali N_{\ph})$ which acts on $H_{\ph}$ by
\bean
 \l(\om)(a_{\ph})=[(\om\circ \k\pd id)\cop(a)]_{\ph},
\eean
and can also be written
\bea
 \l(\om)=(\om\circ\k\pd id)(W).
\lb{pfr}
\eea
$\l$ is called the {\em Fourier representation} of ${\K}$. Its image
$\l(M_*)=\hat M$ is a \vN algebra on which it is true that
\bea
\l(\om)(\om'_{\ph})=(\om\ast \om')_{\ph},
\lb{dpl}
\eea 
where $\om'_{\ph}\in H_{\ph}$ is the unique vector such that
$(\om'_{\ph}|a_{\ph})=\bra a^*,\om'\ket$, for all $a\in\cali N_{\ph}$, and
for every $\om'\in I_{\ph}=\{\om\in M_*\,|\, sup\{|\bra \om,a^*\ket|,\ a\in
M,\ \ph(a^*a)\le 1\}<\infty\}$. Actually, the above condition together with
the definition of $I_{\ph}$ is a generalization of the condition of square
integrability for $\om'\in M_*$.

The dual of the Kac algebra ${\K}$ is based on the image $\hat M$ of the Fourier
representation, $\hat K=(\hat M,\hat\cop,\hat \k,\hat\ph)$. The dual
involution $^o$ goes to $\hat S_{\hat\ph}=\hat
J_{\hat\ph}\hat\cop_{\hat\ph}^{1/2}$, analogously to its dual. The operator
$W$ is unitary and its adjoint is given by (with $\hat J_{\hat\ph}=\hat J$
from now on)
\bea
W^*=(\hat J\pd J)\circ W\circ(\hat J\pd J).
\lb{wd}
\eea
Its dual is then written 
\bea
\hat W=\s\circ W^*\circ\s,
\lb{wh}
\eea
and the dual coproduct is given by the dual version of \rep{cpw}, or under
the form
\bea
 \hat\cop(\om)=\hat W({\1}\pd\om)\hat W^*.
\lb{dwcp}
\eea 

Furthermore, the dual antipode is defined by $\hat \k(\l(\om))=\l(\om\circ
\k)$ or by $\hat \k(\om)=J\om^oJ$, its canonical implementation on $H_{\ph}$.
Dualizing this last formula, we get a new formula for $\k$ in
terms of $\hat J$:
\bea
 \k(a)=\hat Ja^*\hat J.
\lb{rccow}
\eea
The dual weight $\hat\ph$ on $\hat M$ is the normal, faithful semi-finite
weight canonically associated to the left Hilbert algebra $(I_{\ph}\cap
I_{\ph}^o)_{\ph}$ and is given, for $X\in\hat M^+$, by \cite[2.1.1 and
3.5.3.]{ensc}
\bean
\hat{\ph}(X)=\left\{
\begin{array}{cl}
\Vert\om\Vert^2\ &\ \mbox{if there exists }\om\in (I_{\ph}\cap
I_{\ph}^o)_{\ph}\mbox{ such that }X=\hat\p(\om)\\
+\infty\ &\ \mbox{otherwise,}
\end{array}
\right.
\eean
where $\hat\p$ is the canonical representation (acting on the left by the
algebra product) of $(I_{\ph}\cap I_{\ph}^o)_{\ph}$ on $H_{\ph}$. Finally,
the Hilbert space $H_{\hat\ph}$ is identified with $H_{\ph}$. This Kac
algebra obviously has also a predual $\hat M_*$ and a Fourier representation
$\hat\l$.

From \rep{pfr} and by the fact that $\l$ is an involutive representation,
\rep{mein}, it follows by using \rep{wh} that 
\bean
 \l(\om)=(\om\pd id)(W^*)=(id\pd\om)(\hat W).
\eean
If we compare this formula with \rep{gor}, we get $\hat W$ as the generator
of $\l$. Furthermore, if we apply \rep{apr} to $\hat W^*=\s\circ W\circ\s$, we
obtain $(\cop\pd id)(W)=({\1}\pd W)(\s\pd{\1})({\1}\pd W)(\s\pd{\1})$, and from the
relation \rep{canimcop} it follows that $W$ satisfies the {\em pentagonal
relation}
\bean
 ({\1}\pd W)(\s\pd {\1})({\1}\pd W)(\s\pd {\1})(W\pd {\1})=(W\pd {\1})({\1}\pd W).
\eean

The duality mapping between ${\K}$ and $\hat{{\K}}$ is then performed first by
passing from $M$ to its predual $M_*$ and then to $\hat M$ via the Fourier
representation $\l$, which is faithful \cite[Chap.~4]{ensc}. Since $\l$ is
generated by $\hat W$ (and, by duality, $\hat\l$ is generated by $W$), we
understand the essential role played by the operator $W$ in Kac duality. The
Kac duality is then the fact that $\hat{\hat{{\K}}}$ is isomorphic to
${\K}$.

\subsection{The Abelian Kac Algebra of a Group}

Given a separable l.c.\frs group $G$, there are two Kac algebras of special
significance. The first is defined on the \vN algebra $L^{\infty}(G)$ of
(classes of almost everywhere defined) measurable and essentially bounded
functions on $G$ \cite{kad} This means that, for every $f\in L^{\infty}(G)$,
there exists a smallest number $C$ $(0\le C <\infty)$ such that $|f(x)|\le C$
locally almost everywhere. This number $C$ is just the $ess.sup.$ (essential
supremum) of $f$. The norm is given by
\bean
\Vert f\Vert_{\infty}=ess.sup.|f(x)|
\eean
and the involution by $f^*(x)=\ovl{f(x)}$. This algebra acts on the Hilbert
space $L^2(G)$ by pointwise multiplication. This Hilbert space has the $L^2$
scalar product
\bean
 (f|g)=\int_G dx\, f(x)\ovl{g(x)},
\eean
and the norm $\Vert f\Vert_2=\sqrt{(f|f)}$, where $dx$ is the left invariant
Haar measure on $G$. Then $L^{\infty}(G)$, with the operations and conditions
of the following list, is the {\em Abelian Kac algebra} of $G$,
${\K}^a(G)=(L^{\infty}(G),\cop,\k,\ph_a)$:
\bse
 f\cdot g(x)&=&f(x)g(x);\\ 
 {\1}&=&1,\ \mbox{such that}\ 1(x)=1\ \fal\ x\in G;\\ 
 \cop(f)(x,y)&=&f(xy);\\ 
 \k(f)(x)&=&f(\inv x);\\
 \ph_a(f)&=&\int_G dx\ f(x),\hs 1 f\in L^{\infty}(G)^+.
\lb{kag}
\ese
Here $C(G)$ is the algebra of continuous functions with compact support on
$G$, whose product is the convolution (see below). 
The Haar weight is in fact a trace, simply the left invariant Haar measure on
$G$. In consequence, the modular operator is reduced to simple
multiplication by 1, $\cop_{\ph_a}=1$, and the modular group is trivial,
$\s_t^{\ph_a}=id$. The underlying Abelian Hopf-von Neumann algebra is denoted
${\H}^a(G)$.

We also have in this case, for $F\in C(G\pd G)$ and $f\in C(G)$,
\bean
&& WF(x,y)=F(x,xy),\hs 1 W^*F(x,y)=F(x,\inv xy),\\ 
&&\hs 2\hat Jf(x)=\frac 1{\sqrt{\cop_G x}}\ovl{f(\inv x)},
\eean 
where $\cop_G$ is the modular function on $G$ to be defined a few steps below.
From these data we can use relation \rep{dfw} to compute the Fourier
representation for ${\K}^a(G)$,
\bean
 (W(f\pd g)|h\pd l)&=&\int_G dx\,f(x)\ovl{h(x)}\int_G dy\,g(xy)\ovl{l(y)}\\
&=&\int_G dz\,g(z)\int_G dx\,\ovl{h(x)}f(x)\ovl{l(\inv xz)},\\
&=&(g|\l(\om_{hf})l).
\eean
We conclude that $\l(\om_{hf})l(z)=\int_G dx\,h(x)\ovl{f(x)}l(\inv
xz)$ or, taking into account that $\om_{h,f}=h\ovl f\in
L^{\infty}(G)_*$ by \rep{dlf}, that
\bea
 \l(\om_{hf})=\int_G dx\,\om_{hf}(x)L(x).
\lb{forep}
\eea

Since ${\K}^a(G)$ acts on $L^2(G)$, it follows from $\ph_a(f^*f)=\int_G dx\
|f(x)|^2<\infty$ that the GNS construction is given by inclusion, with $\cali
N_{\ph_a}=L^{\infty}(G)\cap L^2(G)$ and $\cali M_{\ph_a}=L^{\infty}(G)\cap
L^1(G)$. The predual is $L^{\infty}(G)_*=L^1(G)$ and, as anticipated,
$I_{\ph_a}=L^1(G)\cap L^2(G)$ is the space of square integrable functions on
the predual $L^1(G)$, on which we now concentrate.

Given a left invariant Haar measure $dx$ on a l.c. group G, the space of
(classes of) functions defined almost everywhere and integrable on G,
$L^1(G,dx)$, is the {\em convolution Banach algebra} of G with as product the
convolution \cite{rei}
\bea
 (f*g)(x)=\int_G dy\ f(y)g(y^{-1}x),
\lb{conv}
\eea
involution
\bea
 f^*(x)=\cop_G\inv x\ovl{f(\inv x)}
\lb{ika}
\eea
and norm $\|| f\||=\int_G dx\ |f(x)|$. $\cop_G:G\apli \Rm$ is a positive
and continuous homomorphism of groups called {\em modular function}:
\bean
\cop_G e&=&1\\ 
\cop_G(xy)&=&\cop_G x\cop_G y.
\eean
If $\m_l$ and $\m_r$ are left- and right-invariant (Haar) measures on G, that
is, $\m_l(xE)=\m_l(E)$, $\m_r(Ey)=\m_r(E)$, for every Borel set $E$, the
function $\cop_G$ relates them by the Radon-Nikod\'ym derivative \cite{rei}
\bea
\frac{d\m_l(x)}{d\m_r(x)}=\cop_G x.
\lb{dfm}
\eea 
When $\cop_G\equiv 1$, the two measures coincide and G is {\em
unimodular}. Changing variables in \rep{conv} and using the identity
\bea
\frac{d\m_l(\inv x)}{d\m_l(x)}=\cop_G\inv x,
\lb{reimem}
\eea
the convolution can also be written in terms of the right invariant measure
(see \rep{dfm}) as $(f\ast g)(x)=\int_G d\m_r(y)\, f(x\inv y)g(y)$.

The algebra $C(G)$ of continuous functions with compact support is dense in
$L^1(G)$, which explains its appearance in some definitions. $G$ is
discrete if and only if $L^1$ has a unit. Otherwise, it has only left and
right approximate units. In general, the algebra $L^1(G)$ is nothing more
than an ideal in the following unital algebra. To every $f\in L^1(G)$ we
associate a measure $d\m(x)$ by $d\m(x)=f(x)dx$. This association implements
an involutive isometry between the Banach algebras $L^1(G)$ and $M^1(G)$, the
unital involutive algebra of all bounded complex Borel measures on G with
convolution given by $(\m *\n)(f)=\int_{G\times G}\ f(xy)d\m(x)d\n(y)=\int_G
\ f(x)d(\m *\n)(x)$, where the unit is the Dirac measure at the identity of
$G$, $\dtd e$. With the notation $\check f(x)=\ovl{f(\inv x)}$, which we
shall be using from now on, the involution is given by
$\m^*(f)=\ovl{\m(\check f)}$ \cite{rei}.

A representation $U$ of $G$ on $H$ is also a representation of $M^1(G)$
and is written 
\bean
 \m\map U(\m)=\int_G\ d\m(x)\ U(x),
\eean
whose restriction to $L^1(G)$ is non-degenerate,
\bea
\lb{lire}
 f\map U(f)=\int_G dx\ f(x) U(x).
\eea
There is, in fact, a bijective correspondence between the unitary irreducible
representations of G and the non-degenerate representations of $L^1(G)$. In
particular, to the left-regular representation of G corresponds the operator
\bea
L(f)\equiv\hat f=\int_G dx f(x)L(x),
\lb{lore}
\eea
whose restriction to $\cali N_{\ph_a}=L^2(G)\cap L^{\infty}(G)$ is just that
derived earlier as the Fourier representation of ${\K}^a(G)$ and denoted
$\l(f)$. Furthermore, when restricted to the space $L^1(G)\cap L^2(G)$,
$L(f)$ acts by convolution:
\bea
 L(f)g=f\ast g,\qquad g\in L^1(G)\cap L^2(G). 
\lb{lorc}
\eea
$L(f)$ also satisfies:
\bean
L(f\ast g)&=&\hat f\cdot\hat g;\\ 
L(f^*)&=&\int_G dx\,\cop_G\inv x\ovl{f(\inv x)}L(x)=\int_G
dx\,\ovl{f(x)}L^{\dagger}(x)\\ &=&\hat f^{\dagger}.
\eean
The operators which constitute the image of $\l$ form the \vN algebra
$\widehat{L^{\infty}(G)}$.

The Abelian \Cal $C_o(G)$ of complex functions vanishing at infinity, with
norm $\|| f\||=sup|f(x)|$ and involution given by the complex conjugation, has
as its dual the algebra $M^1(G)$, the duality pairing being given by
\bean
 \m(f)= \bra \m,f\ket=\int_G d\m(x)f(x)=\int_G dx\ g(x)f(x),
\eean
if $d\m(x)=g(x)dx$. The same duality relation holds between the pair
$L^{\infty}\supset C_o(G)$ and $L^1\subset M^1(G)$ as a linear functional on
the latter, since $L^1(G)^*=L^{\infty}(G)$. While we have
$L^{\infty}(G)_*=L^1(G)$, the dual of $L^{\infty}$ is not $L^1$, but just
contains it \cite{kad}.

\subsection{The Symmetric Kac Algebra of a Group}
\lb{ks}

The \vN algebra $\widehat{L^{\infty}(G)}$, generated by left-regular
representations of a l.c. group G, is denoted $\cali M(G)$. Their generators
$L(x),\ x\in G,$ act on $L^2(G)$ by
\bean
[ L(x)f](y)=f(\inv xy).
\eean
The norm is given by $\||L(x)\||=sup\{\||L(x)f\||,\ \|| f\||_2= 1\}$ and the
involution by hermitian conjugation. The product in $\cali M(G)$ is the
representation of the group product, $L(x)L(y)=L(xy)$, with unit ${\1}=L(e)=I$. 
Every element (operator) in $\cali M(G)$ is a linear combination
of all generators, with functions in $L^1(G)$ as coefficients,
\bea
\lb{ftop}
\hat f=\int_G dx\ f(x)L(x),\hs 1 f\in L^1(G),
\eea
just the image of the Fourier representation of ${\K}^a(G)$ given
in \rep{forep}. These operators act on $L^2(G)$ by
\bea
 [\hat f g](x)=\int_G dy\ f(y) [L(y)g](x)=\int_G dy\ f(y)g(\inv yx).
\lb{fhac} 
\eea
If we further restrict to $g\in L^2(G)\cap L^1(G)$, \rep{fhac} turns out to
be just the convolution \rep{lorc}. The product of operators is written as
\bea
\hat f\cdot\hat g&=&\int_G dx\int_G dy\ f(x)g(y) L(xy)=\int_G dz\int_G dx\
f(x)g(\inv xz) L(z)\non\\ &=&\int_G dz\ (f*g)(z) L(z).
\lb{opp}
\eea

With the operator product \rep{opp}, ${\K}^s(G)=(\cali M(G),\cop,\k,\ph_s)$ is
the {\em symmetric Kac algebra} of the group G. The other operations are
\bse
 &&\cop L(x)=L(x)\pd L(x);\slb{ksaa}\\ 
 &&\k(L(x))=L(\inv x)=\inv L(x)=L^{\dagger}(x);\\ 
 && \ph_s(a)=\left\{
\begin{array}{cl}
 \Vert f\Vert^2_2 & \mbox{if}\ a=\hat f^{\dagger}\cdot\hat f\\
+\infty & \mbox{otherwise.}
\end{array}
\right.\hs 1 a\in\hat M^+ \slb{ksad}
\lb{ksa}
\ese
Just for completeness and better visualization of the structure, we rewrite the
above expressions in terms of the linear combinations \rep{ftop}, whose product
has been given in \rep{opp}:
\bse
 &&\cop(\hat f)=\int_G dx\, f(x)\,L(x)\pd L(x);\\
&&\k(\hat f)=\int_G dx\, f(x)\,L(\inv x)=\int_G dx\,\cop\inv x f(\inv
x)\,L(x);\\ &&\ph_s(\hat f)=f(e),\hs 1 f\in C(G)\ast C(G).\slb{ksd}
\ese

If $F\in C(G\times G)$, we have
\bse
&&\hat WF(x,y)=F(\inv yx,y),\hs 1\hat W^*F(x,y)=F(yx,y),\\
&&\hs 2 Jf(x)=\ovl{f(x)},\ f\in C(G).
\lb{waj}
\ese
In order to see how $\hat W$ generates $\l$, let us consider the space
$L^2(G,L^2(G))$ of $L^2$-valued functions on $G$. It turns out to be
isomorphic to $L^2(G)\pd L^2(G)$ by the association $\f(y)(x)=F(x,y)$,
where $\f(y)\in L^2(G)$. Actually, $\l$ is generated by the left-regular
representation $L$ of $G$, whose action on $\f(y)$ can be written,
with the help of \rep{waj}, as
\bean
 [L(y)\f(y)](x)=\f(y)(\inv yx)&=&F(\inv yx,y)\\
&=&[\hat WF](x,y).
\eean
This shows that the generator $L$, as a bounded function from $G$ to $\cali
B(L^2(G))$, can be seen as the operator $\hat W\in\cali B(L^2(G))\pd
L^{\infty}(G)$.

The modular operator on ${\K}^s(G)$ is given by the Radon-Nikod\'ym derivative
$\cop_{\ph_s}=\frac{d\ph_a}{d(\ph_a\circ\k)}$ of the trace $\ph_a=\m_l$ on
${\K}^a(G)$. Since by a quick calculation one obtains $\m_l\circ\k=\m_r$, it
turns out from the definition of the modular function \rep{dfm} that the
modular operator is just $\cop_G$. The modular function acts on $L^2(G)$ by
pointwise multiplication and, for $f\in H_{\ph_s}$,
\bea
[\s_t^{\ph_s}(L(x))f](y)&=&[\cop_G^{it}L(x)\cop_G^{-it}f](y)\non\\
 &=&(\cop_G y)^{it}[L(x)\cop_G^{-it}f](y)\non\\
 &=&(\cop_G y)^{it}(\cop_G(\inv x y))^{-it}f(\inv xy)\non\\
 &=&(\cop_G x)^{it}f(\inv xy),
\lb{mgks}
\eea
which gives $\s_t^{\ph_s}(L(x))=(\cop_G x)^{it}L(x)$.

As the base space for the dual of ${\K}^s(G)$, the predual $\cali M(G)_*$ is the
\vN algebra of functionals $\hat{\om}_{f,g}:\cali M(G)\apli{\C}$ such that
$\hat{\om}_{f,g}(\hat h)=(\hat h(f)|g)$, which is isomorphic to the {\em
Fourier algebra} $A(G)$ of those functions $h$ which can be written in the form
$h=f\ast\check g,\ f,g\in L^2(G)$. Their identification comes from \rep{dlf}
and is given through the function $\hat{\om}_{f,g}(x)\equiv\bra L(\inv
x),\hat{\om}_{f,g}\ket=(f\ast\check g)(x)$. The Fourier representation in this
case also follows from \rep{dfw} and \rep{waj},
\bean
(\hat W(f\pd g)|h\pd l)&=&\int_G dy\,g(y)\ovl{l(y)}\int_G
dx\,\ovl{h(x)}f(\inv yx)\\ &=&\int_G dy\,g(y)\ovl{(h\ast\check f)(y)l(y)}\\
&=&\int_G dy\,g(y)\ovl{\hat{\om}_{h,f}(y)l(y)},
\eean
from which we get $\hat\l(\hat{\om}_{h,f})l=\hat{\om}_{h,f}l$. It turns out
that the Fourier representation is the identity. By the Cauchy-Schwarz
inequality we obtain that $A(G)$ is contained in $L^{\infty}(G)$ (its normic
closure in this algebra is just $C_o(G)$), and from \rep{dpl} and the last
expression we obtain that it has the usual $L^{\infty}$-pointwise product as
operation.

\subsection{Decomposition into Irreducibles}
\lb{decom}

The reducible representations of a type~I group $G$ can be decomposed into
irreducible representations in a unique way \cite{mack1}. However, the previous
knowledge of the unitary dual $\hat G$ of $G$ is necessary to the actual
realization of the decomposition. The dual is a space consisting of equivalence
classes of unitary irreducible representations with its Mackey-Borel structure
and a Plancherel (this name will be justified below) measure associated to the
Haar measure on $G$ \cite{dix2}. The Plancherel measure and the type~I
Mackey-Borel structure will allow us to sum (or integrate) on $\hat G$ without
ambiguities \cite{mack1,mack}. We proceed to obtain the decomposition in this
section. We will take the regular representation case as a starting point and
arrive at the decomposition of the \vN algebra it generates and of the Hilbert
spaces they act upon (see, for example, \rep{hsd}). For the left-regular
representations, the decomposition can be written in the form
\bean
 L=\int_{\hat G}^{\+} d\m(\al)\, T_{\al},
\eean
where $\al\in\hat G$ and $d\m(\al)$ is a Plancherel measure. This corresponds
to the direct integral decomposition $\cali M(G)=\int^{\+}_{\hat
G}d\m(\al)\,\cali M_{\al}(G)$ of the von Neumann algebra underlying ${\K}^s(G)$.
Since the operators $T_{\al}(x)$ provide irreducible representations of $G$,
there should be an analogous decomposition of the representation of $L^1(G)$,
\bea
 L(f)=\int_{\hat G} d\m(\al)\ T_{\al}(f),
\lb{pfou}
\eea
with each summand given by
\bea
T_{\al}(f)\equiv\hat f_{\al}=\int_{G} dx\ f(x)T_{\al}(x).
\lb{fou}
\eea
This gives a new aspect to formula \rep{ftop},
\bea
L(f)=\int_{\hat G} d\m(\al)\int_{G} dx\ f(x)T_{\al}(x).
\lb{ftopa}
\eea
Formula \rep{fou} is the generalized Fourier transform of $f\in L^1(G)$,
whose outcome is the operator-valued function $\hat f_{\al}$ on $\hat G$. Its
image belongs to the von Neumann algebra $\cali M_{\al}(G)$, which acts on
the Hilbert space $H_{\al}$ such that $L^2(G)=\int_{\hat G}^{\+} d\m(\al)\,
H_{\al}$.

As regards the weight $\ph_s$, it is a trace if and only if $G$, or
${\K}^s(G)$, is unimodular. We can easily show it using \rep{ksd} and
\rep{opp}. Restricting to $f\in L^1(G)\cap L^2(G)$, we get
\bse
&&\ph_s(\hat f\cdot\hat f^{\dagger})=(f\ast f^*)(e)=\int_G dx\,\cop_G
x\,|f(x)|^2\slb{mqp}\\ &&\ph_s(\hat f^{\dagger}\cdot\hat f)=(f^*\ast
f)(e)=\int_G dx\,|f(x)|^2,\slb{mqq}
\lb{mqs}
\ese
where we have also used \rep{reimem} to obtain \rep{mqq}. 
In the unimodular case we have the Plancherel formula, which involves a 
well-defined decomposition for $\ph_s=\Tr$, as explained in
Ref.~\cite{dix2}. In the general case, where symmetric Kac algebras of
a non-unimodular type~I group are considered, we can suppose also the weight
$\ph_s$ to be decomposed according to
\bea
 \ph_s=\int_{\hat G}^{\+} d\m(\al)\,\ph_{\al}.
\lb{hwdec}
\eea 
It will be sometimes useful to extend abusively the domain of $\ph_s$ to the
generators $L(x)$, which can be regarded as left-regular representations of
the Dirac measures $\d_x\in M^1(G):\ L(x)=\int_G dy\, \d_x(y)L(y)$. In
this sense we write
\bea
\d_e(x)=\ph_s(L(x))=\int_{\hat G} d\m(\al)\,\ph_{\al}(T_{\al}(x)),
\lb{gdm}
\eea
which is to be regarded as the explicit general expression for the Dirac delta
distribution on the group.

Going further, from \rep{mqp} we can write, for $f\in L^1(G)\cap L^2(G)$, two
expressions: on one hand, $\ph_s(\hat f^{\dagger}\cdot\hat f)=\int_G
dx\,|f(x)|^2$; on the other hand, $\ph_s(\hat f^{\dagger}\cdot\hat
f)=\int_{\hat G} d\m(\al)\,\ph_{\al}[(\hat f^{\dagger}\cdot\hat f)_{\al}]$. We
obtain, consequently,
\bean
\int_G dx\,|f(x)|^2=\int_{\hat G} d\m(\al)\,\ph_{\al}[(\hat
f^{\dagger}\cdot\hat f)_{\al}]
\eean
as a generalization of the Plancherel formula, where
\bean
(\hat f^{\dagger}\cdot\hat f)_{\al}=\hat f_{\al}^{\dagger}\cdot\hat
f_{\al}=\int_G dx\,(f^*\ast f)(x)T_{\al}(x).
\eean
Since $f\in L^2(G)$, it follows that, for {\em almost all} $\al$,
$\ph_{\al}[\hat f^{\dagger}_{\al}\cdot\hat f_{\al}]<\infty$, and we can
conclude that $\hat f_{\al}\in\cali N_{\ph_{\al}}$ for {\em almost all}
$\al$. Here and in the following {\em almost all} $\al$ includes that set of
representations whose complement in the unitary dual has zero Plancherel
measure, that is, the support of this measure. It is
generally identified with the set of higher dimensional representations. For
example, in the half-plane canonical group case they are just the
infinite-dimensional $T_{\pm}$. From \rep{hwdec} we also have that, for
almost all $\al$, the $\ph_{\al}$ are normal, faithful and semi-finite
weights on $\cali M_{\al}(G)$.

With the above weight decomposition we are able to write out the inverse of
the generalized Fourier transform \rep{fou}. From \rep{ksd} it follows that
\bean
f(x)=\ph_s[L^{\dagger}(x)\hat f],
\eean
whose decomposition
\bea
 f(x)=\int_{\hat G}d\m(\al)\,\ph_{\al}[T_{\al}^{\dagger}(x)\hat f_{\al}]
\lb{gfti}
\eea 
gives $f(x)$ in terms of the operator-valued function $\hat f_{\al}$ on $\hat
G$. Writing
\bea
 f_{\al}(x)\equiv\ph_{\al}[T_{\al}^{\dagger}(x)\hat f_{\al}]
\lb{dgfti}
\eea 
and recalling that $f\in L^1(G)$, it follows from \rep{gfti} that $\int_G
dx\,|f_{\al}(x)|<\infty$ for almost all $\al$, that is, that $f_{\al}\in
L^1(G)$ for almost all $\al$.

Notice that the generalized Fourier transform defined in \rep{fou} and
\rep{gfti} is faithful as a map between $G$ and its dual if and only if the
Plancherel measure accounts for every element of $\hat G$. Since the
Plancherel measure is concentrated on the highest dimensional
representations, it may happen in some cases, like those of the Heisenberg
\cite{fol,tay} and the special half-plane \cite{gur} groups, that there are
irreducible inequivalent representations on $\hat G$ which are missed in
formulas \rep{pfou} and \rep{gfti}.

We have up to now collected the decompositions of $\cali M(G)$, of its
generators $L(x)$, of the Hilbert space $L^2(G)$ and of the Haar weight
$\ph_s$. The question coming naturally to the mind is whether or not these
(irreducible) components constitute a Kac algebra. The answer is negative,
because the components $\ph_{\al}$ of $\ph_s$ are not Haar weights in the
Hopf-von Neumann algebra ${\H}_{\al}(G)$ generated by $T_{\al}(x)$. The structure
of ${\H}_{\al}(G)$ for fixed $\al\in\hat G$ comes straight from the decomposition
of $L$:
\bse 
&&T_{\al}(x)T_{\al}(y)=T_{\al}(xy)\\
&&I=T_{\al}(e)\\
&&\cop_{\al}(T_{\al}(x))=T_{\al}(x)\pd T_{\al}(x)\slb{hcop}\\
&&\k_{\al}(T_{\al}(x))=T_{\al}^{\dagger}(x).
\lb{hvna}
\ese
If we take $\ph_{\al}$ and try to verify the second axiom for Haar weights,
for example, we get from the two sides of \rep{syma}:
\bean
(id\pd\ph_{\al})[(I\pd
T_{\al}^{\dagger}(y))\cop_{\al}(T_{\al}(x))]&=&\ph_{\al}(T_{\al}(\inv
yx))\,T_{\al}(x),\\
\k_{\al}(id\pd\ph_{\al})[\cop_{\al}(T_{\al}^{\dagger}(y))(I\pd
T_{\al}(x))]&=&\ph_{\al}(T_{\al}(\inv yx))\,T_{\al}(y),
\eean
which implies that axiom if and only if $x=y$. Since there is no
warrant that $\ph_{\al}(T_{\al}(\inv yx))$ would have as outcome $x=y$ (we
have instead \rep{gdm}), we conclude that $\ph_{\al}$ is not Haar.
Conversely, by the Haar weight axioms, it can be proved that a weight $\ph'$
is a Haar weight on ${\H}_{\al}(G)$ if and only if $\ph'(T_{\al}(x))=\d_e(x)$.

The elements of ${\H}_{\al}(G)$ are written 
\bea
 \t_{\al}(f)\equiv\hat f_{\al}=\int_G dx\, f(x)\,T_{\al}(x),\hs 1 f\in L^1(G).
\lb{wefo}
\eea 
This means that they are the images of the inequivalent irreducible
representations of $L^1(G)$ corresponding to the $T_{\al}$ representations of
$G$. In analogy with the relation between the representation $L$ and $\l$,
where the latter is a restriction of $L$ to $\cali
N_{\ph_a}=L^{\infty}(G)\cap L^2(G)$, formula \rep{wefo} is regarded as a
restriction of formula \rep{fou} to the respective decomposition of $\cali
N_{\ph_a}$, that is, as an $\al$-component $\t_{\al}$ of the Fourier
representation $\l$. It is thus natural to look for its generator. We shall,
in what follows, restrict ourselves to separable and type~I semi-direct
product groups of the type $G=S\pds N$, where $N$ is an Abelian normal
subgroup and $S$ is a unimodular group. This restriction on the group will
provide more explicit formulas for the Kac algebra decomposition, while
retaining enough generality to allow the consideration not only of the
half-plane example envisaged here but also of other cases of physical
interest, like the Euclidean motion group $E(2)$, the correct canonical group
for the phase space of the circle \cite{ish}. Elements of $G$ will be
denoted by $x=(s,n),\ y=(r,l)$, etc., the identity by $(e,0)$ and the product
by $(s,n)(r,l)=(sr,n+\f_s(l))$, where $\f_s$ is a homomorphism on $N$, the
action of $S$. In this case, a generalization of what was presented in
section~\ref{irrur} regarding the group $\RmR$ is provided by Mackey's theory
of induced representations applied to semi-direct product groups (see also
Refs.~\cite{sugi,gur,bar}). If $V_y$ are irreducible representations
(characters) of $N$, labeled by $y\in\hat N$, that theory says that the
Hilbert space $H_y$ is formed by those functions $f_y$ which satisfy:
\bit
\item the map $(s,n)\in G\map f_y(s,n)\in{\C}$ is measurable;
\item $f_y((s,n)(e,l))=V_y^{-1}(l)\,f_y(s,n)$;
\item $\int_{G/N} d\m(s)\, |f_y(s,n)|^2 <\infty$,
\eit
where $d\m(s)$ is a G-invariant measure on $G/N\sim S$ (notice that $H_y$
differs from $H_{\al}$ in that the label $\al$ represents classes of
inequivalent representations while $y$ represents all (irreducible)
representations). The action of $(s,n)\in G$ on $r\in S$, denoted, $(s,n)\cdot
r$, is defined by taking the $S$-component of the product
$(s,n)(r,0)=(sr,n)=(sr,0)(e,\f_{sr}^{-1}(n))$, according to the decomposition
$(s,n)=(s,0)(e,\f_s^{-1}(n))$ of $G$ in terms of its parts $S$ and $N$, that
is, $(s,n)\cdot r\equiv sr\in S$. By the same decomposition, the second
condition implies that the functions $f_y$ can be written as
\bea
f_y(s,n)=V_y^{-1}(\f_s^{-1}(n))\x(s),\hs 1 f_y(s,0)\equiv\x(s)\in L^2(S).
\lb{fdeai}
\eea
Actually, \rep{fdeai} expresses an isomorphism between $H_y\ (H_{\al})$ and
$L^2(S)$ for each label $y\,(\al)$. In what follows we will suppose that the
irreducible representations have been already classified, that is, we will work
on $H_{\al}$. On these spaces the irreducible unitary induced representations
$T_{\al}$ of $G$ by $V_{\al}$ are given by
\bea
 [T_{\al}(x)\x](s)=V_{\al}(\s(\inv x\cdot s;x))\x(\inv x\cdot s),
\lb{irepi}
\eea
where $\s(r;x)$ is a ``gaugefied'' cocycle on $G$, or a $(S,G)$-cocycle
relative to the invariant class of $d\m(s)$ \cite{var,gur}, that is, a Borel
map $\s:S\times G\apli N$ which satisfies
\bea
&& \s(r;e)=0;\\
&& \s(y\cdot r;x) - \s(r;xy) + \s(r;y)=0.
\lb{idqdocs}
\eea
It is given explicitly by $\s(r;s,n)=\f_{sr}^{-1}(n)$. ``Gaugefied'' cocycles
appear already in usual Quantum Mechanics, even in its discretized version, in
which the Euclidean phase space is replaced by ${{\Z}}_n\pd{{\Z}}_n$ \cite{alga}.

Formula \rep{wefo} will have an important role in our work. It generalizes
Weyl's formula \cite{wey} in the sense that it associates ($L^1$) functions on
the group to irreducible operators on the subgroup $S$. In order to find
explicitly the generator of the representation $T_{\al}(\t_{\al})$ of $L^1(G)$,
we consider functions $\ps\in L^2(G,H_{\al})\sim L^2(G,L^2(S))$ and, putting
$\ps(x)(s)=G(s,x),\ x\in G$ we define an isomorphism between the
$L^2(S)$-valued functions on $G$ and the space $L^2(S)\pd L^2(G)$. Now, the
induced irreducible representations on $\ps(x)\in L^2(S)$ are given by
\bea
[T_{\al}(x)\ps(x)](s)&=&V_{\al}(\s(\inv x\cdot s;x))\ps(x)(\inv x\cdot s)\non\\
&=&V_{\al}(\s(\inv x\cdot s;x))G(\inv x\cdot s,x)
\equiv [\hat W^{\al}G](s,x).
\lb{goir}
\eea
This shows that the generator of $\t_{\al}$, the function $T_{\al}$ in
$L^{\infty}(G,\cali B(L^2(S)))$, can be seen as an operator $\hat
W^{\al}\in\cali B(L^2(S))\pd L^{\infty}(G)$. Operator $\hat W^{\al}$ is
analogous to $\hat W$ not only in what regards the generation of Fourier
representations but also because it implements the coproduct
\rep{hcop}. This can be shown by recalling the definitions of the induced
representations $T_{\al}$ on $H_{\al}(G)$ and $L^2(S)$, and comparing the two
actions below,
\bean
&& [\hat W^{\al}(I\pd T_{\al}(z))\hat W^{\al *}G_{\al}](s,x)=\\ &&\hs 2=
V_{\al}^{-1}(-\s(\inv x\cdot s;x))V_{\al}^{-1}(\s(\inv x\cdot s;\inv
zx))\,G(\inv z\cdot s,\inv zx)\\ &&\cop_{\al}
T_{\al}(z)G_{\al}(s,x)=V_{\al}^{-1}(\s(s;\inv x))\,G(\inv z\cdot s,\inv zx).
\eean
The left hand sides above turn out to be equal if we substitute in the right
hand sides the identity $\s(\inv x\cdot s;\inv zx)=\s(s;\inv z)+\s(\inv
x\cdot s;x)$, straightforwardly obtained from \rep{idqdocs}. We have thus
that $W^{\al}$ is the fundamental operator of ${\H}_{\al}(G)$. It is easily
verified that it satisfies the pentagonal relation.

Turning back to \rep{wefo}, we obtain by \rep{irepi} that the operators $\hat
f_{\al}$ act on $L^2(S)$ by
\bea
[\hat f_{\al}\x](r)=\int_G d^l\m(s,n)\, V_{\al}(\s(\inv sr;s,n))f(s,n)\x(\inv
sr).
\lb{acof}
\eea
Since the right invariant measure on $G$ is the product of the invariant
measures on $S$ and $N$, we have $d^l\m(s,n)=\cop(s,n)\, d\m(s)\,d\m(n)$.
After the change of variables $\inv sr=t$ and by Fubini's theorem, \rep{acof}
reads
\bean
 [\hat f_{\al}\x](r)=\int_S d\m(t)\, K_f^{\al}(r,t) \x(t)
\eean
in terms of the kernel $K_f^{\al}(r,t)$ given by
\bean
K_f^{\al}(r,t)=\int_N d\m(n)\,\cop(r\inv t,n)\,V_{\al}(\s(t;r\inv
t,n))f(r\inv t,n).
\eean
Introducing a kernel will enable us to write out an explicit formula for the
weight $\ph_{\al}$. Also the following result will help: {\em the modular
function of a semi-direct product group $G=S\pds N$ is only a function on
$S$}. This is proved by using $d^r\m(\inv xy)=\cop_G x d^r\m(y)$ and the
invariance of the measures on $S$ and $N$:
\bean
 d^r\m((s,n)^{-1}(r,l))&=&d^r\m(\inv sr,\f_s^{-1}(l-n))\\
&=&d\m(\inv sr)\pex d\m(\f_s^{-1}(l-n))\\
&=&\left|\frac{\partial\f_s^{-1}(l)}{\partial l}\right|\,d\m(r)\pex d\m(l)\\
&=&\cop_G(s,n)\,d^r\m(r,l),
\eean
that is, 
\bea 
\cop_G(s,n)=\left|\frac{\partial\f_s^{-1}(l)}{\partial l}\right|\equiv\cop(s),
\lb{mfin}
\eea
and, in particular, $\cop_G(e,n)=\cop(e)=1$. In $\G$ we have
$d^r(\l,a)=\frac{d\l}{\l} da$ and the left-invariant measure is easily
verified to be $d^l(\l,a)=d\l da$, which implies that
$\cop_{\G}(\l,a)=\cop(\l)=\l$. Turning back to the general case, a trace can
be introduced on ${\H}_{\al}(G)$ by
\bse
\Tr_{\al}(\hat f_{\al})&=&\int_S d\m(t)\, K_f^{\al}(t,t)\slb{trdef}\\
&=&\int_G d\m(t)\,d\m(n)\,V_{\al}(\s(t;e,n))f(e,n).
\ese
Formula \rep{trdef} is a good trace definition because the kernels
satisfy
\bean
 \int_S d\m(t)\, K_f^{\al}(r,t)K_g^{\al}(t,s)=K_{f\ast g}^{\al}(r,s),
\eean 
which implies $\Tr_{\al}(\hat f_{\al}^*\hat f_{\al})=\Tr_{\al}(\hat
f_{\al}\hat f_{\al}^*)$. 
We will now introduce an explicit decomposition of the Haar weight in terms of
the trace:
\bea
\ph_{\al}(\hat f_{\al})&\equiv&\Tr_{\al}(\cop\hat f_{\al})=\int_S
d\m(t)\,\cop(t)K_f^{\al}(t,t)\non\\ &=&\int_G
d\m(t)\,d\m(n)\,\cop(t)\,V_{\al}(\s(t;e,n))f(e,n)\non\\ &=&\int_G
d^l\m(t,n)\,V_{\al}(\s(t;e,n))f(e,n),
\lb{wafo}
\eea
where $\cop$ is given by \rep{mfin}. Clearly it is not a trace. For example,
in the half-plane group we have
\bean
 \ph_{\pm}(\hat f_{\pm})=\frac 1{2\p}\int_{\G} d\l da\,e^{\pm ia\l}f(1,a),
\eean
which is a decomposition of the Haar weight $\ph_s$, since
\bean
\ph_s(\hat f)=\sum_{\pm}\ph_{\pm}(\hat f_{\pm})=\frac
1{\p}\int_{\Rm\times{\R}}d\l da\,\cos(a\l)f(1,a)=f(1,0).
\eean
Computing $\ph_{\al}(T_{\al}(r,l))$, which should be $\d_e(r,l)$ if $\ph_{\al}$
were a Haar weight, the formula above provides another way to see why that does
not happen:
\bean
\ph_{\al}(T_{\al}(r,l))=\d_e(r)\int_S d\m(t)\cop(t)\,V_{\al}(\s(t;e,l)).
\eean 
For $\G$ this gives
\bean
\ph_{\pm}(T_{\pm}(\l,a))=\frac{\d_1(\l)}{2\p}\int_{\Rm}d\r\, e^{\pm
ia\r}=\frac{\d_1(\l)}{2\p}\left(\pi\d(a)\pm\frac ia\right),
\eean
while $\sum_{\pm}\ph_{\pm}(T_{\pm}(\l,a))=\d_1(\l)\d(a)$.
Furthermore, from formula \rep{wafo} the function \rep{dgfti} reads explicitly
\bean
f_{\al}(r,l)=\int_G d^l\m(t,n)\,V_{\al}(\s(t;e,n))f((r,l)(e,n)),\hs 1f\in L^1(G).
\eean

We turn now our attention to the predual of $\cali M_{\al}(G)$ which, in
analogy with the Fourier algebra, we call $A_{\al}(G)$. As in that case,
we introduce the representative function $\hat\om^{\al}_{\x,\ch}(x)$ of this
new algebra by
\bean
\hat\om^{\al}_{\ch,\x}(x)\equiv\bra\hat\om^{\al}_{\ch,\x},T_{\al}^{\dagger}(x)
\ket,
\eean
which, by definition of $\hat\om^{\al}_{\x,\ch}$ as a linear form, is given by
the scalar product
\bea
\hat\om^{\al}_{\ch,\x}(s,n)&=&(T_{\al}^{\dagger}(s,n)\ch|\x)_{L^2(S)}=
(\ch|T_{\al}(s,n)\x)_{L^2(S)}\non\\
 &=&\int_S d\m(t)\, V_{\al}^{-1}(\s(\inv st;s,n))\ch(t)\check{\x}(\inv ts),
\lb{omah} 
\eea
where, we recall, $\check{\x}(s)=\ovl{\x(\inv s)}$. The product on $A_{\al}(G)$
is obtained by the duality shown in \rep{dpcp} from the coproduct on
${\H}_{\al}(G)$ and is the same Abelian pointwise product of $A(G)$, since
\rep{hcop} is symmetric, and of the same kind of the coproduct on $K^s(G)$.
Another point, the involution on $A_{\al}(G)$ is straightforwardly seen from
\rep{rei} to be just the complex conjugation. These facts show that
$A_{\al}(G)$ is very similar to $A(G)$, differing from it only in that their
elements should be written according to \rep{omah} and depend on the labels
$\al\in\hat G$. As $A(G)$ is contained in $L^{\infty}(G)$, this suggests that
$A_{\al}(G)$ be contained in some space alike. To see it better, we compute the
modulus of $\hat\om^{\al}_{\ch,\x}(s,n)$ and, using the Cauchy-Schwarz
inequality, obtain
\bean
 |\hat\om^{\al}_{\ch,\x}(s,n)|&=&|(\ch|T_{\al}(s,n)\x)|\\
 &\le&\Vert\ch\Vert_2\Vert T_{\al}(s,n)\x\Vert_2\\
 &=&\Vert\ch\Vert_2\Vert\x\Vert_2<\infty,
\eean
since $T_{\al}$ is unitary and $\ch,\x\in L^2(S)$. Thus,
$|\hat\om^{\al}_{\ch,\x}(s,n)|$ is essentially bounded and we can say
that $A_{\al}(G)$ is contained in $L^{\infty}(G)$. The index $\al$ just
indicates the dependence on the labels in $\hat G$.

Using the explicit form of the generator $\hat W^{\al}$ we can determine the
representation $\hat\t_{\al}$ of $\cali M_{\al}(G)_*=A_{\al}(G)$ by formula
\rep{dfw},
\bea
(\hat W^{\al}(\x,f)|\ch\pd g)_{L^2(S)\pd
L^2(G)}=(f|\hat\t_{\al}(\hat\om^{\al}_{\ch,\x})g)_{L^2(G)},
\lb{mfbwl}
\eea
where $\x,\ch\in L^2(S)$ and $f,g\in L^2(G)$. The left hand side gives
\bean
(\hat W^{\al}(\x,f)|\ch\pd g)&=&
\int_G d\m^l(s,n)\, f(s,n)\int_S d\m(t)\,V_{\al}(\s(\inv
st;s,n))\ovl{\ch(t)}\x(\inv st)\,\ovl{g(s,n)}\\ &=&\int_G d\m^l(s,n)\,
f(s,n)\,\ovl{\hat\om^{\al}_{\ch,\x}(s,n)g(s,n)}.
\eean
Comparison with the right hand side of \rep{mfbwl} yields $\hat\t_{\al}=id$.
Taking $\hat\t_{\al}$ as an $\al$-component of $\hat\l$ and recalling that the
dual of $K^s(G)$ is built on the image of $\hat\l$, we conclude that the dual
of ${\H}_{\al}(G)$ is built on $\hat\t_{\al}(A_{\al}(G))\subset
L^{\infty}(G)$, that is, the dual is contained in ${\H}^a(G)$. We also obtain
from \rep{dpl} that, while $\cali M_{\al}(G)$ acts on $L^2(S)$, $A_{\al}(G)$
acts on $L^2(G)$ (by the pointwise product), which explains the asymmetry of
the double scalar product in \rep{mfbwl}. The dual version of that formula,
\bea
(W^{\al}(f,\x)|g\pd\ch)_{L^2(G)\pd
L^2(S)}=(\x|\t_{\al}(\om_{g,f})\ch)_{L^2(S)},
\lb{muel}
\eea 
which is also asymmetric, involves the representation \rep{wefo} and the
operator $W^{\al}=\s\circ\hat W^{\al *}\circ\s$ ($f\in L^2(G),\ \x\in L^2(S)$),
\bean
[W^{\al}(f,\x)](s,n;r)=V_{\al}^{-1}(\s(r;s,n))\,f(s,n)\x(sr).
\eean
The left hand side of \rep{muel} then gives 
\bean
(W^{\al}(f,\x)|g\pd\ch)&=&
\int_S d\m(t)\,\x(t)\int_G d\m^l(s,n)\,V_{\al}^{-1}(\s(\inv
st;s,n))\ovl{g(s,n)}f(s,n)\,\ovl{\ch(\inv st)}\\ &=&\int_S
d\m(t)\,\x(t)\ovl{\int_G d\m^l(s,n)\,\om_{g,f}(s,n)[T_{\al}(s,n)\ch](t)},
\eean
where the first equality involves a change of variables in $S$ and we have
identified $\om_{g,f}=g\ovl f$ by \rep{dlf}. Comparison of this result with the
right hand side of \rep{muel} corroborates formula \rep{wefo} for $\t_{\al}$.
If we also introduce in $A_{\al}(G)$ the $\al$-component of the coproduct in
${\H}^a(G)$, this will be implemented by $W^{\al}$, in the same way that $\hat
W^{\al}$ implements \rep{hcop}.

\section{Quantization on the Half-Plane}\lb{qohpl}

We have now at hand a powerful structure to describe quantization. A
generalization of the Weyl-Wigner correspondence prescription is incorporated
in Kac (group) duality. Our objective in this section is to specialize to the
half-plane case the last section results, particularly those concerning the
decomposition of the Kac algebras. We shall show that the Hopf-von Neumann
algebra generated by the irreducible operators, together with its dual, does
provide the framework in which quantization can be described as an
irreducible component in Kac duality theory.

Starting from the group $\G$ as the closest algebraic entity associated to
the half-plane phase space, we necessarily have to consider -- if we are
thinking about duality -- the two Kac algebras ${\K}^s(\G)$ and ${\K}^a(\G)$. The
decomposition of the first according to the dual space $\widehat{\G}$ leads
to the family of Hopf-von Neumann algebras ${\H}_{\al}(\G)$, which inherit most
of their structure from the Kac algebras they come from. Though group duality
is lost at the Hopf-von Neumann level, a well-defined formula for the
decomposition of the Fourier representation persists. Adaptation of formula
\rep{wefo} and its inverse \rep{gfti}, although not representing a bijection
between the group and its dual, provides a well-defined mapping between
functions and irreducible operators.

The cohomological differences between the complete-plane and the half-plane
cases come from the necessity of central-extending the special canonical
group of the former to the Heisenberg group in order to provide a faithful
momentum map between its associated Lie and Poisson algebras \cite{ish}.
Despite these differences, we can regard the Weyl formula as being formula 
\rep{wefo} for a fixed value of the label $\al$ in terms of
the Planck constant, $\al=\al(h)$. This is easily confirmed if we recall that
the unitary dual of the Heisenberg group $H(3)$ is {\em almost equal} (in the
{\em almost everywhere} sense) to any of the $\Om$-projective unitary dual
\cite{mack2} of the bidimensional translation group ${\R}^2$. By the Stone-von
Neumann theorem \cite{tay}, the former is equal to $({\Z}-\{0\})\cup{\R}^2$, while
the latter is just ${\Z}-\{0\}$. We have shown in a separate paper \cite{ralas2}
that the Weyl-Wigner formalism can be described in terms of duality of {\em
projective Kac algebras}. In such a projective duality framework for ${\R}^2$,
Weyl's formula comes from an expression analogous to \rep{wefo} for the
decomposition of the respective Fourier representation, namely
\bea
\hat f_{\n}=\int_{{\R}^2}dx dy\, f(x,y)\, e^{-i\n(y\hat q+x\hat p)},\hs
1\n\in{\Z}-\{0\}\lb{weyl},
\eea
where $\hat q,\ \hat p$ are the usual position and momentum operators. 
Comparing \rep{weyl} with Weyl's formula, we get immediately $\n=\hbar^{-1}$.
Actually, the only formal difference between \rep{wefo} and the original
Weyl's formula, or \rep{weyl}, is that the latter is written in terms of
unitary irreducible {\em projective} operators instead of the linear ones which
appear in formula \rep{wefo}. This is a consequence of the necessity of a
central extension in the complete-plane case. That is, Quantum Mechanics is a
theory on a particular Hilbert space and its operators generate a particular
Hopf-von Neumann algebra whose label in the Kac algebra decomposition is just
a point in the support of the Plancherel measure on the unitary dual
space of the group involved. In the half-plane, as observed at the end of
section~\ref{canon}, there is no need for a central extension, since the
cohomology group $H^2(\G,{\R})$ is trivial and consequently projective and
linear representations are cohomologically equivalent. This enables us to use
the simpler, linear representations. Thus, the analogue of Weyl's
correspondence formula for the half-plane group is given by \rep{wefo} for a
fixed value of the labels $\pm$. From the dual
$\widehat{\G}=\{+\}\cup\{-\}\cup{\R}$ we have that
\rep{wefo} is in this case given by
\bse
&&\hat f_{\pm}=\int_{\G}d\l da\, f(\l,a)\,T_{\pm}(\l,a)\slb{ioba} \\
&&f(y)=\int_{\G} d\l da\, f(\l,a)\,\ch_y(\l),\hs 1 y\in{\R}.\slb{iobb}
\ese
Notice that the Hopf-von Neumann algebras ${\H}_{\pm}(\G)$ of operators
\rep{ioba} are quite different from those of functions \rep{iobb}
and denoted ${\H}_y(\G)$. They are Abelian for each $y$ and their direct
integral over ${\R}$ can be identified with ${\H}^a({\R})$.

Recalling that the labels $\pm$ correspond to an uncountable infinity of
equivalent representations in the support of the Plancherel measure, and taking
into account the physical dimensions of the elements of $\G$ ($[\l]=$ length,
$[a]=$ momentum, $[\hbar]=[a\l]=[\hat\p]=$ action), we take $\pm\hbar^{-1}$ for
the representatives of each class instead of $\pm$, and fix the value of the
label to be $+\hbar^{-1}$. The quantization map is then
given by (we will write $\hbar$ instead of $+\hbar^{-1}$ when it appears only
as an indicative of class)
\bea
\hat f_{\hbar}=\int_{\G}d\l da\, f(\l,a)\, e^{\frac i{\hbar}a\hat\r}e^{-\frac
i{\hbar}\ln(\l/\l_o)\hat\p},
\lb{qtmp}
\eea
where $\l_o$ is a constant with dimension of length. The self-adjoint operators $\hat\r$ and $\hat\p$ act on the subspace of $L^2(\Rm)$-functions vanishing at $0$ and $\infty$  by
\bean
\hat\r\x(\r)&=&\r\x(\r)\\
\hat\p\x(\r)&=&-i\hbar\r\frac{\partial\x(\r)}{\partial\r},
\eean
and satisfy the commutation relation
\bean
 [\hat\r,\hat\p]=i\hbar\hat\r.
\eean
The function $f(\l,a)$ is recovered from the operator $\hat f_{\hbar}$ by the
inverse mapping \rep{gfti}
\bea
f(\l,a)=\sum_{\pm\hbar^{-1}}\ph_{\pm\hbar}[T_{\pm\hbar}^{\dagger}(\l,a)
\hat f_{\pm\hbar}]=\sum_{\pm\hbar^{-1}}f_{\pm\hbar}(\l,a),
\lb{qtin}
\eea
where
\bea
 \ph_{\pm\hbar}[\hat f_{\pm\hbar}]=\frac 1{2\p\hbar}\int_{\G}d\r db\,
e^{\pm\frac i{\hbar}b\r}f(1,b),
\lb{wohvn}
\eea
which gives 
\bea
 f_{\pm\hbar}(\l,a)=\frac 1{2\p\hbar}\int_{\G}d\r db\,
e^{\pm\frac i{\hbar}b\r}f((\l,a)(1,b)).
\lb{qtind}
\eea
Eq.~\rep{qtin} explicits the fact that the classical $L^1$-function $f$ has
contributions from {\em almost all} irreducible representations, while
Eq.~\rep{qtind} is just the projection of that function into one of its
``components''. 

The symmetric but non-Abelian Hopf-von Neumann algebra ${\H}_{\hbar}(\G)$
generated by the irreducible operators $T_{\hbar}(\l,a)$ is then the operator
algebra of Quantum Mechanics on the half-plane. Its trivial structure is
analogue to that given in \rep{hvna}. On this algebra there is also defined a
weight given by the plus sign in \rep{wohvn}, which is an irreducible component
of the Haar weight on ${\K}^s(\G)$, as shown in the previous section.

On the dual Abelian Hopf-von Neumann algebra ${\H}^a(\G)$, a typical
$A_{\hbar}(\G)$-function is that given by
\bean
\hat\om^{\hbar}_{\ch,\x}(\l,a)&=&(\ch |T_{\hbar}(\l,a)\x)_{L^2(\Rm)}\\
&=&\iRm{\r}\, e^{-\frac i{\hbar}a\r}\ch(\r)\ovl{\x(\r/\l)}.
\eean
If we put $\ch=\x$, $\hat\om^{\hbar}_{\x,\x}\equiv W_{\x}^{\hbar}$ is to be
interpreted as a generalization of the Wigner distribution function
for the half-plane associated to the state $\x$. This is justified, for if we
compute the expectation value of the operator $\hat f_{\hbar}$ in the state
$\x$, it is given by
\bean
\bra\hat f_{\hbar}\ket_{\x}=(\x|\hat f_{\hbar}\x)=\int_{\G} d\l
da\,f(\l,a)\,W_{\x}^{\hbar}(\l,a).
\eean
This makes clear the role of $W_{\x}^{\hbar}$ as a {\em quantum probability
density}, the same role played by the Wigner distribution in the Euclidean
phase space. But notice that things here are quite different from the
complete-plane case and this function does not share most of the properties
the usual Wigner distribution is known to satisfy. The differences are due to
a lack of connection between functions in $A_{\hbar}(\G)$, like
$W_{\x}^{\hbar}$, and $L^1$-functions, or the respective operator in $\cali
M_{\hbar}(\G)$. Banach duality is not able to provide an explicit
correspondence between these two spaces when the group is not Abelian
self-dual like $\G$. In the complete-plane case, ${\R}^2$ is self-dual and the
Banach duality turns out to be just the double Fourier transform on the phase
space. Furthermore, the Abelian algebras $L^1$ and $L^{\infty}$ over ${\R}^2$
are isomorphic by the uniqueness of the Fourier transform \cite{rei}. This
gives rise to the well-known formulas of the Weyl correspondence
\cite{fol,tay,hil}, of which the Wigner distribution function is a particular
case corresponding to the density operator. And, since there is no need to
consider projective representations in the half-plane case, no 2-cocycle
arises, that is, neither the convolution algebra $L^1(\G)$ is twisted nor the
pointwise product algebra $L^{\infty}(\G)$ is deformed by any kind of
Moyal-like product.

\section{Final Comments}

The Weyl-Wigner prescription is based on the Pontryagin duality of Fourier
transformations. It calls attention to the central role of duality in
quantization, though it only can be expected to hold in the particular case of
Abelian canonical groups. We have been concerned with the impact that general
Fourier duality can have in quantization. The stage-sets for Fourier analysis
are neither groups nor homogeneous spaces, but Kac algebras. General Fourier
duality requires actually a pair of algebras and we have considered such a pair
of ``symmetric'' and ``Abelian'' Kac algebras for a particular, type~I but
non-Abelian and non-unimodular, canonical group. The decomposition of the first
has led to some Hopf-von Neumann algebras on the group, which we have
recognized as the natural algebraic arenas where duality plays its role and,
furthermore, where we can find out how far it is possible to go with the
Weyl-Wigner approach as a guideline to quantize general systems.

The open half-plane which we have examined is perhaps the simplest case
presenting some novel, deep features. It is still globally Euclidean -- though
no more a vector space. Although we have been restricted to a case in which the
special canonical group and the phase space manifolds coincide, the group
non-triviality requires new algebraic structures. In particular, the group
involved being non-unimodular, it is no more a trace which is at work, but its
generalization allowing noncommutative integration -- a weight. To connect Kac
duality and Weyl quantization, we must restrict ourselves to a specific
irreducible representation of the group involved. The operators in that
representation generate a Hopf-von Neumann algebra which participates in the
irreducible decomposition of the symmetric Kac algebra of the canonical group.
It is in general impossible to obtain an explicit correspondence between the
$L^1$- and $L^{\infty}$-functions on it. From the point of view of Fourier
duality, this is standard -- Fourier transforms in general map $L^1$-functions
into $\hat G$-valued operators ($L^{\infty}$-functions only if $G$ is Abelian).
In the Wigner formalism this corresponds to a failure in the correspondence
between the Wigner distribution and its density. Generalized Wigner
distribution functions only make sense if related to density operators and, as
such, they are defined as the expectation values of the irreducible operators
$T_{\hbar}$. Summing up, generalized Fourier duality in the case treated here
does provide a prescription for the quantization of $L^1$-functions on phase
space via a generalized Weyl's formula. Although it is possible to recover the
quantizable c-number function from the correspondent Weyl operator, the
correspondence is not at all complete, since we cannot relate it to its dual
$L^{\infty}$-function.

The conclusion is that general Fourier duality does provide a firm guide to
quantization, though imposing severe restrictions to the simple-minded
expectations born from the results concerning those very simple systems for
which the phase space is a vector space. Since this duality in only achieved in
the Kac algebra framework, we also conclude that, for quantization purposes,
algebraic structures beyond groups must be considered.

\section*{Acknowledgments}

The authors would like to thank Prof. M. Enock for useful comments on Kac
algebras intricacies. They also would like to thank CAPES and CNPq for
financial support.

\appendix
\section*{Coadjoint Orbits of $\G$}

For the sake of completeness we shall show here that the half-plane is the
unique non-trivial homogeneous symplectic manifold by the action of the group
$\G=\RmR$. To do that, we realize the group as a $2\times 2$ matrix group by
the correspondence
\bean
(\l,a)\map\left(
\begin{array}{cc}
\l^{-1/2}\ &\ a\l^{1/2}\\
	0\ &\ \l^{1/2}
\end{array}
\right).
\eean
The Lie algebra $\lG$ can be accordingly realized if we define its
generators by
\bean
L=\frac d{dt}(e^t,0)|_{t=0}=\left(
\begin{array}{cc}
-1/2\ &\ 0\\ 
0\ &\ 1/2
\end{array}
\right)\hs 1
A=\frac d{dt}(1,t)|_{t=0}=\left(
\begin{array}{cc}
0\ &\ 1\\ 
0\ &\ 0
\end{array}
\right).
\eean
A quick computation is enough to verify that $L, A$ realize the algebra $\lG$,
\bean
 [A,L]=A.
\eean 

To get the adjoint action of $\G$ on $\lG$, we write an arbitrary element
$X$ of the algebra as $X=X^AA+X^LL$, and compute
\bean
Ad_{(\l,a)}X=(\l,a)X(\l,a)^{-1}=(aX^L+\l^{-1}X^A)A + X^LL.
\eean
Now, to obtain the coadjoint action of $\G$ on $\lG^*$, we first find a dual
basis to $X_{\nu}=\{A,L\}$ in $\lG^*$ through the duality pairing
$\bra\th^{\mu},X_{\nu}\ket=\Tr(\th^{\mu}X_{\nu})=\d^{\mu}_{\nu}$, and get
\bean
\th^L=\left(
\begin{array}{cc}
-1\ &\ 0\\ 
 0\ &\ 1
\end{array}
\right)\qquad
\th^A=\left(
\begin{array}{cc}
0\ &\ 0\\ 
1\ &\ 0
\end{array}
\right).
\eean
Computing the coadjoint action on an element
$\et=\et_L\th^L+\et_A\th^A\in\lG^*,\et_{\nu}\in{\R}$, we find
\bean
 [Ad^*_{(\l,a)}\et](X)&\equiv&\bra\et,Ad^{-1}_{(\l,a)}X\ket\\
&=&\bra\et,X^LL + \l(X^A-aX^L)A\ket\\ &=&X^L\et_L + \l(X^A-aX^L)\et_A.
\eean 
To obtain it for any $X\in\lG$, we compare with
\bean
[Ad^*_{(\l,a)}\et](X)=(Ad^*_{(\l,a)}\et)_LX^L+(Ad^*_{(\l,a)}\et)_AX^A,
\eean
which gives finally
\bea
Ad^*_{(\l,a)}\et=\l\et_A\th^A+(\et_L-a\l\et_A)\th^L.
\lb{aca}
\eea

The orbits of this action on $\lG^*$ are given for all $(\l,a)\in\G$.
Analyzing the coefficients of $\th^A$ and $\th^L$ in \rep{aca}, we conclude
that there are basically two kind of orbits: those for which $\et_A\neq
0$; and those for which $\et_A=0$. In the first case the coefficient of
$\th^A$ is never zero but that of $\th^L$ can assume any value in ${\R}$.
This characterizes two half-planes, one for $\et_A>0$ and the other for
$\et_A<0$. In the second case, $(\et_A=0)$, we have
$Ad^*_{(\l,a)}\et=\et_L\th^L$, which means that these orbits consist of the
infinity of isolated points in the line $\et_A=0$. This concludes our
analysis, showing that $\G$ has two 2-dimensional orbits diffeomorphic to the
half-plane $({\lG^*}_{\pm})$ and an uncountable infinity of
0-dimensional orbits $({\lG^*}_0)$ in $\lG^*$.

To conclude, we can also compute the Kirillov symplectic form on the orbits
passing by $\et$ by the formula $\om_{\et}=\frac
12C_{\mu\nu}^{\s}\et_{\s}\th^{\mu}\pex\th^{\nu}$. Since $C_{AL}^A=1$, we
have on ${\lG^*}_{\pm}\sim\Rm\times{\R}$
\bean
\om_{\pm}=\et_A\th^A\pex\th^L.
\eean
The symplectic form $\om$ used in section~\ref{canon} is obtained from $\om_-$
above through the realization
\bean
A\map \del p&\qquad& L\map p\del p-x\del x\\
\th^A\map dp + pd\ln x&\qquad& \th^L\map - d\ln x,
\eean
and with $\et_A=-x,\ x\in\Rm$.

\end{document}